\documentclass[11pt]{article}
\usepackage{amssymb,amsthm,amsmath}
\def\qed{\hfill $\Box$}

\newcommand\pf{\smallbreak\noindent \texttt{Proof}. }

\begin{document}

\newtheorem{thm}{Theorem}[section]
\newtheorem{prop}[thm]{Proposition}
\newtheorem{lem}[thm]{Lemma}
\newtheorem{cor}[thm]{Corollary}
\newtheorem{ex}[thm]{Example}
\renewcommand{\thefootnote}{*}

\title{\bf On the inner structure of 3-dimensional Leibniz algebras}

\author{\textbf{L.A.~Kurdachenko, O.O.~Pypka}\\
Oles Honchar Dnipro National University, Dnipro, Ukraine\\
{\small e-mail: lkurdachenko@gmail.com, sasha.pypka@gmail.com}\\
\textbf{I.Ya.~Subbotin}\\
National University, Los Angeles, USA\\
{\small e-mail: isubboti@nu.edu}}
\date{}

\maketitle

\begin{abstract}
Let $L$ be an algebra over a field $F$ with the binary operations $+$ and $[,]$. Then $L$ is called a \textit{left Leibniz algebra} if $[[a,b],c]=[a,[b,c]]-[b,[a,c]]$ for all $a,b,c\in L$. We describe the inner structure of left Leibniz algebras having dimension 3.
\end{abstract}



\thispagestyle{empty}

\section{Introduction}
Let $L$ be an algebra over a field $F$ with the binary operations $+$ and $[,]$. Then, $L$ is called a \textit{Leibniz algebra} (more precisely, a \textit{left Leibniz algebra}) if, for all elements $a,b,c\in L$, it satisfies the Leibniz identity:
$$[[a,b],c]=[a,[b,c]]-[b,[a,c]].$$
We will also use another form of this identity:
$$[a,[b,c]]=[[a,b],c]+[b,[a,c]].$$

Leibniz algebras first appeared in the paper of A.~Blokh~\cite{BA1965} where they were called \textit{$D$-algebras}. However, at that time, these works were not in demand. Only after two decades was there a real rise in interest towards Leibniz algebras. It happened thanks to the rediscovery of these algebras by J.-L.~Loday~\cite{LJ1993} (see also~\cite[Section~10.6]{LJ1992}), who used the term ``Leibniz algebras'' since it was Leibniz who discovered and proved the Leibniz rule for differentiation of the product of functions. The main motivation for the introduction of Leibniz algebras was the study of periodicity phenomena in algebraic K-theory. The Leibniz algebras appeared to be naturally related to several areas such as differential geometry, homological algebra, classical algebraic topology, algebraic K-theory, loop spaces, non-commutative geometry, physics and so on. Nowadays, the theory of Leibniz algebras is one of actively developing areas of modern algebra.

Note that Lie algebras are a partial case of Leibniz algebras. Conversely, if $L$ is a Leibniz algebra in which $[a,a]=0$ for every element $a\in L$, then it is a Lie algebra. Thus, Leibniz algebras can be seen as a non-commutative generalization of Lie algebras.

The theory of Leibniz algebras has been intensely developing in many different directions. Some results of this theory were presented in the recent book~\cite{AOR2020}.

One of the first steps in the theory of Leibniz algebras is the description of algebras with small dimensions. Unlike Lie algebras, the situation with Leibniz algebras of dimension 3 is very diverse. Leibniz algebras of dimension 3 are mostly described. The description of Leibniz algebras of dimensions 4 and 5 is quite complex. The list of papers devoted to these studies is quite large and we will not give it here in full. We only note that the Section~3.1 of book~\cite{AOR2020} is devoted to study of right Leibniz algebras having dimension 3. The investigation of Leibniz algebras having dimension 3 was carried out in articles~\cite{CILL2012,DMS2014,KKO2010,RRM2018,RR2012,YaV2019}. Some of these papers use the language of the right Leibniz algebras, whilst others use the language of left Leibniz algebras. Basically, the description is reduced to determining the structural constants of these algebras. However, the structural constants do not always give an idea of the internal structure of these algebras. Elucidation of the structure requires some additional analysis. The overall picture seems fragmented. Thus, the articles dealt with Leibniz algebras over concrete fields of real, complex, $p$-adic numbers, etc. Therefore, we cannot decide if the internal structure of these algebras really contains a complete description. For example, when passing from the field of rational numbers to the field of real numbers, some types of algebras disappear altogether. Furthermore, some sections describe right Leibniz algebras, while others describe left Leibniz algebras. Moreover, only structural constants were found. None of the articles considered the internal structure. So, under these circumstances, in order to observe the entire scope, there are two options that arise: analyze these articles and really make sure everything is done there or do it yourself all over again. The second option is preferable since it is better to complete the list within your own framework if it turns out that not all types of algebras have been considered. More importantly, we are interested in the description of the inner structure, not just structural constants. Of course, you can extract some information about the structure from structural constants, but it is more logical to obtain the structural constants in the process of describing the inner structure.

Therefore, in the current article, we present a description of left Leibniz algebras having dimension 3, focusing on clarifying their structure and obtaining structural constants by passing. This consideration of the structure of Leibniz algebras of dimension 3 is carried out over an arbitrary field $F$, and when studying these specific types of algebras, additional natural restrictions on the field $F$ appear. These restrictions are very significant in some cases. Some types of algebras can exist only if sufficiently strict restrictions are imposed. Our goal was the most detailed description of these algebras, reflecting all the nuances of their structure.

\section{Main results.}
Let $L$ be a Leibniz algebra over a field $F$. Then $L$ is called \textit{abelian} if $[a,b]=0$ for every elements $a,b\in L$. In particular, an abelian Leibniz algebra is a Lie algebra.

If $A,B$ are subspaces of $L$, then $[A,B]$ will denote a subspace generated by all elements $[a,b]$ where $a\in A$, $b\in B$. As usual, a subspace $A$ of $L$ is called a \textit{subalgebra} of $L$, if $[a,b]\in A$ for every $a,b\in A$. It follows that $[A,A]\leqslant A$. A subalgebra $A$ of $L$ is called a \textit{left} (respectively \textit{right}) \textit{ideal} of $L$, if $[b,a]\in A$ (respectively $[a,b]\in A$) for every $a\in A$, $b\in L$. In other words, if $A$ is a left (respectively right) ideal, then $[L,A]\leqslant A$ (respectively $[A,L]\leqslant A$). A subalgebra $A$ of $L$ is called an \textit{ideal} of $L$ (more precisely, \textit{two-sided ideal}) if it is both a left ideal and a right ideal.

Every Leibniz algebra $L$ possesses the following specific ideal. Denote by $\mathrm{\mathbf{Leib}}(L)$ the subspace generated by the elements $[a,a]$, $a\in L$. It is not hard to prove that $\mathrm{\mathbf{Leib}}(L)$ is an ideal of $L$. The ideal $\mathrm{\mathbf{Leib}}(L)$ is called the \textit{Leibniz kernel} of algebra $L$.

We note the following important property of the Leibniz kernel:
$$[[a,a],x]=[a,[a,x]]-[a,[a,x]]=0.$$

The \textit{left} (respectively \textit{right}) \textit{center} $\zeta^{\mathrm{left}}(L)$ (respectively $\zeta^{\mathrm{right}}(L)$) of a Leibniz algebra $L$ is defined by the rule:
$$\zeta^{\mathrm{left}}(L)=\{x\in L|\ [x,y]=0\ \mbox{for each element }y\in L\}$$
(respectively,
$$\zeta^{\mathrm{right}}(L)=\{x\in L|\ [y,x]=0\ \mbox{for each element }y\in L\}).$$
It is not hard to prove that the left center of $L$ is an ideal, but that is not true for the right center. Moreover, $\mathrm{\mathbf{Leib}}(L)\leqslant\zeta^{\mathrm{left}}(L)$ so that $L/\zeta^{\mathrm{left}}(L)$ is a Lie algebra. The right center is a subalgebra of $L$ and, in general, the left and right centers are distinct (see, for example,~\cite{KOP2016}).

The \textit{center} $\zeta(L)$ of $L$ is defined by the rule:
$$\zeta(L)=\{x\in L|\ [x,y]=0=[y,x]\ \mbox{for each element }y\in L\}.$$
The center is an ideal of $L$.

Now we define the \textit{upper central series}
$$\langle0\rangle=\zeta_{0}(L)\leqslant\zeta_{1}(L)\leqslant\ldots\zeta_{\alpha}(L)\leqslant\zeta_{\alpha+1}(L)\leqslant\ldots\zeta_{\eta}(L)$$
of a Leibniz algebra $L$ by the following rule: $\zeta_{1}(L)=\zeta(L)$ is the center of $L$, and recursively, $\zeta_{\alpha+1}(L)/\zeta_{\alpha}(L)=\zeta(L/\zeta_{\alpha}(L))$ for all ordinals $\alpha$, and $\zeta_{\lambda}(L)=\bigcup_{\mu<\lambda}\zeta_{\mu}(L)$ for the limit ordinals $\lambda$. By definition, each term of this series is an ideal of $L$.

Define the \textit{lower central series} of $L$
$$L=\gamma_{1}(L)\geqslant\gamma_{2}(L)\geqslant\ldots\gamma_{\alpha}(L)\geqslant\gamma_{\alpha+1}(L)\geqslant\ldots\gamma_{\tau}(L)$$
by the rule: $\gamma_{1}(L)=L$, $\gamma_{2}(L)=[L,L]$, $\gamma_{\alpha+1}(L)=[L,\gamma_{\alpha}(L)]$ for all ordinals $\alpha$ and $\gamma_{\lambda}(L)=\bigcap_{\mu<\lambda}\gamma_{\mu}(L)$ for the limit ordinals $\lambda$.

As usual, we say that a Leibniz algebra $L$ is \textit{nilpotent}, if there exists a positive integer $k$ such that $\gamma_{k}(L)=\langle0\rangle$. More precisely, $L$ is said to be \textit{nilpotent of nilpotency class $c$} if $\gamma_{c+1}(L)=\langle0\rangle$, but $\gamma_{c}(L)\neq\langle0\rangle$. We denote the nilpotency class of $L$ by $\mathrm{\mathbf{ncl}}(L)$.

Define the \textit{lower derived series} of $L$
$$L=\delta_{0}(L)\geqslant\delta_{1}(L)\geqslant\ldots\delta_{\alpha}(L)\geqslant\delta_{\alpha+1}(L)\geqslant\ldots\delta_{\nu}(L)$$
by the following rule: $\delta_{0}(L)=L$, $\delta_{1}(L)=[L,L]$, and recursively $\delta_{\alpha+1}(L)=[\delta_{\alpha}(L),\delta_{\alpha}(L)]$ for all ordinals $\alpha$ and $\delta_{\lambda}(L)=\bigcap_{\mu<\lambda}\delta_{\mu}(L)$ for the limit ordinals $\lambda$. If $\delta_{n}(L)=\langle0\rangle$ for some positive integer $n$, then we say that $L$ is a \textit{soluble} Leibniz algebra.

As usual, we say that a Leibniz algebra $L$ is \textit{finite dimensional} if the dimension of $L$ as a vector space over $F$ is finite.

If $\mathrm{\mathbf{dim}}_{F}(L)=1$, then $L$ is abelian.

If $\mathrm{\mathbf{dim}}_{F}(L)=2$, then we obtain the following types of Leibniz algebras:
\begin{equation*}
\begin{split}
\mathrm{\mathbf{Lei}}_{1}(2,F)&=Fa\oplus Fb\ \mbox{ where }[a,a]=b,[a,b]=[b,a]=[b,b]=0;\\
\mathrm{\mathbf{Lei}}_{2}(2,F)&=Fc\oplus Fd\ \mbox{ where }[c,c]=[c,d]=d,[d,c]=[d,d]=0
\end{split}
\end{equation*}
(see, for example,~\cite{KKPS2017}).

Moving on to Leibniz algebras of dimension 3, we immediately note that we will consider Leibniz algebras, which are not Lie algebras. This means that their Leibniz kernel is non-zero. Then, the factor-algebra over Leibniz kernel has dimension, at most, 2. Note that Lie algebras having dimension at most 2 are soluble. Thus, we obtain the following

\begin{prop}
Let $L$ be a Leibniz algebra over a field $F$. Suppose that $L$ is not a Lie algebra. If $L$ has dimension $3$, then $L$ is soluble.
\end{prop}

For the Leibniz kernel $\mathrm{\mathbf{Leib}}(L)$ of a Leibniz but not a Lie algebra, $L$ having dimension 3 will give us only two possibilities: $\mathrm{\mathbf{dim}}_{F}(\mathrm{\mathbf{Leib}}(L))=1$ and $\mathrm{\mathbf{dim}}_{F}(\mathrm{\mathbf{Leib}}(L))=2$.

First, we will consider the situation when $\mathrm{\mathbf{dim}}_{F}(\mathrm{\mathbf{Leib}}(L))=1$. Immediately, we obtain the following two subcases:

(IA) the center of $L$ includes $\mathrm{\mathbf{Leib}}(L)$;

(IB) the Leibniz kernel of $L$ is not central.

\noindent For each of these subcases, we have the following two possibilities:

(IA1) the factor-algebra $L/\mathrm{\mathbf{Leib}}(L)$ is abelian;

(IA2) the factor-algebra $L/\mathrm{\mathbf{Leib}}(L)$ is not abelian;

\noindent and

(IB1) the factor-algebra $L/\mathrm{\mathbf{Leib}}(L)$ is abelian;

(IB2) the factor-algebra $L/\mathrm{\mathbf{Leib}}(L)$ is not abelian.

Consider these cases.

\begin{thm}\label{T1}
Let $L$ be a Leibniz algebra over a field $F$ having dimension $3$. Suppose that $L$ is not a Lie algebra. If the center of $L$ includes the Leibniz kernel, $\mathrm{\mathbf{dim}}_{F}(\mathrm{\mathbf{Leib}}(L))=1$ and the factor-algebra $L/\mathrm{\mathbf{Leib}}(L)$ is abelian, then $L$ is an algebra of one of the following types.
\begin{enumerate}
\item[\upshape(i)] $\mathrm{\mathbf{Lei}}_{3}(3,F)=L_{3}$ is a direct sum of two ideals $A=Fa_{1}\oplus Fa_{3}$ and $B=Fa_{2}$. Moreover, $A$ is a nilpotent cyclic Leibniz algebra of dimension $2$, $[A,B]=[B,A]=\langle0\rangle$, so that $L_{3}=Fa_{1}\oplus Fa_{2}\oplus Fa_{3}$ where $[a_{1},a_{1}]=a_{3}$, $[a_{1},a_{2}]=[a_{1},a_{3}]=[a_{2},a_{1}]=[a_{2},a_{2}]=[a_{2},a_{3}]=[a_{3},a_{1}]=[a_{3},a_{2}]=[a_{3},a_{3}]=0$, $\mathrm{\mathbf{Leib}}(L_{3})=[L_{3},L_{3}]=Fa_{3}$, $\zeta^{\mathrm{left}}(L_{3})=\zeta(L_{3})=\zeta^{\mathrm{right}}(L_{3})=Fa_{2}\oplus Fa_{3}$, $L_{3}$ is nilpotent and $\mathrm{\mathbf{ncl}}(L_{3})=2$.
\item[\upshape(ii)] $\mathrm{\mathbf{Lei}}_{4}(3,F)=L_{4}$ is a direct sum of ideal $A=Fa_{1}\oplus Fa_{3}$ and subalgebra $B=Fa_{2}$. Moreover, $A$ is a nilpotent cyclic Leibniz algebra of dimension $2$, $[A,B]=Fa_{3}$, $[B,A]=\langle0\rangle$, so that $L_{4}=Fa_{1}\oplus Fa_{2}\oplus Fa_{3}$ where $[a_{1},a_{1}]=[a_{1},a_{2}]=a_{3}$, $[a_{1},a_{3}]=[a_{2},a_{1}]=[a_{2},a_{2}]=[a_{2},a_{3}]=[a_{3},a_{1}]=[a_{3},a_{2}]=[a_{3},a_{3}]=0$, $\mathrm{\mathbf{Leib}}(L_{4})=[L_{4},L_{4}]=\zeta^{\mathrm{right}}(L_{4})=\zeta(L_{4})=Fa_{3}$, $\zeta^{\mathrm{left}}(L_{4})=Fa_{2}\oplus Fa_{3}$, $L_{4}$ is nilpotent and $\mathrm{\mathbf{ncl}}(L_{4})=2$.
\item[\upshape(iii)] $\mathrm{\mathbf{Lei}}_{5}(3,F)=L_{5}$ is a direct sum of ideal $A=Fa_{1}\oplus Fa_{3}$ and subalgebra $B=Fa_{2}$. Moreover, $A$ is a nilpotent cyclic Leibniz algebra of dimension $2$, $[A,B]=\langle0\rangle$, $[B,A]=Fa_{3}$, so that $L_{5}=Fa_{1}\oplus Fa_{2}\oplus Fa_{3}$ where $[a_{1},a_{1}]=[a_{2},a_{1}]=a_{3}$, $[a_{1},a_{2}]=[a_{1},a_{3}]=[a_{2},a_{2}]=[a_{2},a_{3}]=[a_{3},a_{1}]=[a_{3},a_{2}]=[a_{3},a_{3}]=0$, $\mathrm{\mathbf{Leib}}(L_{5})=[L_{5},L_{5}]=\zeta^{\mathrm{left}}(L_{5})=\zeta(L_{5})=Fa_{3}$, $\zeta^{\mathrm{right}}(L_{5})=Fa_{2}\oplus Fa_{3}$, $L_{5}$ is nilpotent and $\mathrm{\mathbf{ncl}}(L_{5})=2$.
\item[\upshape(iv)] $\mathrm{\mathbf{Lei}}_{6}(3,F)=L_{6}$ is a direct sum of ideal $A=Fa_{1}\oplus Fa_{3}$ and subalgebra $B=Fa_{2}$. Moreover, $A$ is a nilpotent cyclic Leibniz algebra of dimension $2$, $[A,B]=[B,A]=Fa_{3}$, so that $L_{6}=Fa_{1}\oplus Fa_{2}\oplus Fa_{3}$ where $[a_{1},a_{1}]=[a_{2},a_{1}]=a_{3}$, $[a_{1},a_{2}]=\alpha a_{3}$ $(\alpha\neq0)$, $[a_{1},a_{3}]=[a_{2},a_{2}]=[a_{2},a_{3}]=[a_{3},a_{1}]=[a_{3},a_{2}]=[a_{3},a_{3}]=0$, $\mathrm{\mathbf{Leib}}(L_{6})=[L_{6},L_{6}]=\zeta^{\mathrm{right}}(L_{6})=\zeta^{\mathrm{left}}(L_{6})=\zeta(L_{6})=Fa_{3}$, $L_{6}$ is nilpotent and $\mathrm{\mathbf{ncl}}(L_{6})=2$.
\item[\upshape(v)] $\mathrm{\mathbf{Lei}}_{7}(3,F)=L_{7}$ is a sum of two nilpotent cyclic ideals $A=Fa_{1}\oplus Fa_{3}$ and $C=Fa_{2}\oplus Fa_{3}$, $[A,C]=[C,A]=\langle0\rangle$, so that $L_{7}=Fa_{1}\oplus Fa_{2}\oplus Fa_{3}$ where $[a_{1},a_{1}]=a_{3}$, $[a_{2},a_{2}]=\beta a_{3}$ $(\beta\neq0)$, $[a_{1},a_{2}]=[a_{1},a_{3}]=[a_{2},a_{1}]=[a_{2},a_{3}]=[a_{3},a_{1}]=[a_{3},a_{2}]=[a_{3},a_{3}]=0$, $\mathrm{\mathbf{Leib}}(L_{7})=[L_{7},L_{7}]=\zeta^{\mathrm{right}}(L_{7})=\zeta^{\mathrm{left}}(L_{7})=\zeta(L_{7})=Fa_{3}$. Moreover, polynomial $X^{2}+\beta$ has no root in field $F$, $L_{7}$ is nilpotent and $\mathrm{\mathbf{ncl}}(L_{7})=2$.
\item[\upshape(vi)] $\mathrm{\mathbf{Lei}}_{8}(3,F)=L_{8}$ is a sum of two nilpotent cyclic ideals $A=Fa_{1}\oplus Fa_{3}$ and $C=Fa_{2}\oplus Fa_{3}$, $[A,C]=Fa_{3}$, $[C,A]=\langle0\rangle$, so that $L_{8}=Fa_{1}\oplus Fa_{2}\oplus Fa_{3}$, $[a_{1},a_{1}]=a_{3}$, $[a_{1},a_{2}]=\alpha a_{3}$, $[a_{2},a_{2}]=\beta a_{3}$ $(\alpha,\beta\neq0)$, $[a_{1},a_{3}]=[a_{2},a_{1}]=[a_{2},a_{3}]=[a_{3},a_{1}]=[a_{3},a_{2}]=[a_{3},a_{3}]=0$, $\mathrm{\mathbf{Leib}}(L_{8})=[L_{8},L_{8}]=\zeta^{\mathrm{right}}(L_{8})=\zeta^{\mathrm{left}}(L_{8})=\zeta(L_{8})=Fa_{3}$. Moreover, polynomial $X^{2}+\alpha X+\beta$ has no root in field $F$, $L_{8}$ is nilpotent and $\mathrm{\mathbf{ncl}}(L_{8})=2$.
\item[\upshape(vii)] $\mathrm{\mathbf{Lei}}_{9}(3,F)=L_{9}$ is a sum of two nilpotent cyclic ideals $A=Fa_{1}\oplus Fa_{3}$ and $B=Fa_{2}\oplus Fa_{3}$ such that $[A,B]=\langle0\rangle$, $[B,A]=Fa_{3}$, so that $L_{9}=Fa_{1}\oplus Fa_{2}\oplus Fa_{3}$ where $[a_{1},a_{1}]=a_{3}$, $[a_{2},a_{1}]=a_{3}$, $[a_{2},a_{2}]=\sigma a_{3}\ (\sigma\neq0)$, $[a_{1},a_{2}]=[a_{1},a_{3}]=[a_{2},a_{3}]=[a_{3},a_{1}]=[a_{3},a_{2}]=[a_{3},a_{3}]=0$, $\mathrm{\mathbf{Leib}}(L_{9})=[L_{9},L_{9}]=\zeta^{\mathrm{right}}(L_{9})=\zeta(L_{9})=\zeta^{\mathrm{left}}(L_{9})=Fa_{3}$. Moreover, polynomial $X^{2}+X+\sigma$ has no root in field $F$, $L_{9}$ is nilpotent and $\mathrm{\mathbf{ncl}}(L_{9})=2$.
\item[\upshape(viii)] $\mathrm{\mathbf{Lei}}_{10}(3,F)=L_{10}$ is a sum of two nilpotent cyclic ideals $A=Fa_{1}\oplus Fa_{3}$ and $B=Fa_{2}\oplus Fa_{3}$ such that $[A,B]=[B,A]=Fa_{3}$, so that $L_{10}=Fa_{1}\oplus Fa_{2}\oplus Fa_{3}$ where $[a_{1},a_{1}]=a_{3}$, $[a_{2},a_{1}]=a_{3}$, $[a_{1},a_{2}]=\tau a_{3}\ (\tau\neq0)$, $[a_{2},a_{2}]=\sigma a_{3}\ (\sigma\neq0)$, $[a_{1},a_{3}]=[a_{2},a_{3}]=[a_{3},a_{1}]=[a_{3},a_{2}]=[a_{3},a_{3}]=0$, $\mathrm{\mathbf{Leib}}(L_{10})=[L_{10},L_{10}]=\zeta^{\mathrm{right}}(L_{10})=\zeta^{\mathrm{left}}(L_{10})=\zeta(L_{10})=Fa_{3}$. Moreover, polynomial $X^{2}+(\tau+1)X+\sigma$ has no root in field $F$, $L_{10}$ is nilpotent and $\mathrm{\mathbf{ncl}}(L_{10})=2$.
\end{enumerate}
\end{thm}
\pf We note that the center $\zeta(L)$ has dimension at most 2. Suppose first that $\mathrm{\mathbf{dim}}_{F}(\zeta(L))=2$. Since $L$ is not a Lie algebra, there is an element $a_{1}$ such that $[a_{1},a_{1}]=a_{3}\neq0$. We note that $a_{3}\in\zeta(L)$. It follows that $[a_{1},a_{3}]=[a_{3},a_{1}]=[a_{3},a_{3}]=0$. Being an abelian algebra of dimension 2, $\zeta(L)$ has a direct decomposition $\zeta(L)=Fa_{2}\oplus Fa_{3}$ for some element $a_{2}$. Put $A=Fa_{1}\oplus Fa_{3}$ and $B=Fa_{2}$, then $B\leqslant\zeta(L)$, so that $B$ is an ideal of $L$. Clearly, $L=Fa_{1}\oplus Fa_{2}\oplus Fa_{3}=A\oplus B$ and $A$ is also an ideal of $L$. Moreover, $A$ is a nilpotent cyclic Leibniz algebra of dimension 2. Thus, we come to the following type of nilpotent Leibniz algebras:
\begin{gather*}
L_{3}=Fa_{1}\oplus Fa_{2}\oplus Fa_{3}\ \mbox{where }[a_{1},a_{1}]=a_{3}, [a_{1},a_{2}]=\\ [a_{1},a_{3}]=[a_{2},a_{1}]=[a_{2},a_{2}]=[a_{2},a_{3}]=[a_{3},a_{1}]=[a_{3},a_{2}]=[a_{3},a_{3}]=0.
\end{gather*}
Note also that $\mathrm{\mathbf{Leib}}(L_{3})=[L_{3},L_{3}]=Fa_{3}$, $\zeta^{\mathrm{left}}(L_{3})=\zeta^{\mathrm{right}}(L_{3})=\zeta(L_{3})=Fa_{2}\oplus Fa_{3}$, $\mathrm{\mathbf{ncl}}(L_{3})=2$.

Suppose now that the center of $L$ has dimension 1. In this case, $\zeta(L)=\mathrm{\mathbf{Leib}}(L)$. Since $L$ is not a Lie algebra, there is an element $a_{1}$ such that $[a_{1},a_{1}]=a_{3}\neq0$. We note that $a_{3}\in\zeta(L)$. It follows that $[a_{1},a_{3}]=[a_{3},a_{1}]=[a_{3},a_{3}]=0$. Then $\zeta(L)=Fa_{3}$. Since $L/\mathrm{\mathbf{Leib}}(L)$ is abelian, for every element $x\in L$ we have $[a_{1},x],[x,a_{1}]\in\zeta(L)\leqslant\langle a_{1}\rangle=Fa_{1}\oplus Fa_{3}$. It follows that subalgebra $\langle a_{1}\rangle$ is an ideal of $L$. Since $\mathrm{\mathbf{dim}}_{F}(\langle a_{1}\rangle)=2$, $\langle a_{1}\rangle\neq L$.

Suppose first that there exists an element $b$ such that $b\not\in\langle a_{1}\rangle$ and $[b,b]=0$. We have $[b,a_{1}]=\gamma a_{3}$ for some $\gamma\in F$. The following two cases appear here: $\gamma=0$ and $\gamma\neq0$. Let $\gamma=0$. Then $[a_{1},b]=\alpha a_{3}$ for some $\alpha\in F$. If we suppose that $\alpha=0$, then $b\in\zeta(L)$. But in this case, $\mathrm{\mathbf{dim}}_{F}(\zeta(L))=2$, and we obtain a contradiction, which shows that $\alpha\neq0$. Put $a_{2}=\alpha^{-1}b$, then $[a_{2},a_{2}]=[a_{2},a_{1}]=0$, $[a_{1},a_{2}]=a_{3}$, and we come to the following nilpotent Leibniz algebra:
\begin{gather*}
L_{4}=Fa_{1}\oplus Fa_{2}\oplus Fa_{3}\ \mbox{where }[a_{1},a_{1}]=[a_{1},a_{2}]=a_{3},\\
[a_{1},a_{3}]=[a_{2},a_{1}]=[a_{2},a_{2}]=[a_{2},a_{3}]=[a_{3},a_{1}]=[a_{3},a_{2}]=[a_{3},a_{3}]=0.
\end{gather*}
Note also that $\mathrm{\mathbf{Leib}}(L_{4})=[L_{4},L_{4}]=Fa_{3}$, $\zeta^{\mathrm{left}}(L_{4})=Fa_{2}\oplus Fa_{3}$, $\zeta(L_{4})=\zeta^{\mathrm{right}}(L_{4})=Fa_{3}$, $\mathrm{\mathbf{ncl}}(L_{4})=2$.

Let $\gamma\neq0$. Put $a_{2}=\gamma^{-1}b$, then $[a_{2},a_{2}]=0$, $[a_{2},a_{1}]=a_{3}$. We have $[a_{1},a_{2}]=\alpha a_{3}$ for some element $\alpha\in F$. If $\alpha=0$, then $a_{2}\in\zeta^{\mathrm{right}}(L)$, and we come to the following nilpotent Leibniz algebra:
\begin{gather*}
L_{5}=Fa_{1}\oplus Fa_{2}\oplus Fa_{3}\ \mbox{where }[a_{1},a_{1}]=[a_{2},a_{1}]=a_{3},\\
[a_{1},a_{2}]=[a_{1},a_{3}]=[a_{2},a_{2}]=[a_{2},a_{3}]=[a_{3},a_{1}]=[a_{3},a_{2}]=[a_{3},a_{3}]=0.
\end{gather*}
Note also that $\mathrm{\mathbf{Leib}}(L_{5})=[L_{5},L_{5}]=Fa_{3}$, $\zeta^{\mathrm{right}}(L_{5})=Fa_{2}\oplus Fa_{3}$, $\zeta(L_{5})=\zeta^{\mathrm{left}}(L_{5})=Fa_{3}$, $\mathrm{\mathbf{ncl}}(L_{5})=2$.

Suppose that $\alpha\neq0$. Then, we come to the following type of nilpotent Leibniz algebras:
\begin{gather*}
L_{6}=Fa_{1}\oplus Fa_{2}\oplus Fa_{3}\ \mbox{where }[a_{1},a_{1}]=[a_{2},a_{1}]=a_{3},\\
[a_{1},a_{2}]=\alpha a_{3}\ (\alpha\neq0),\\
[a_{1},a_{3}]=[a_{2},a_{2}]=[a_{2},a_{3}]=[a_{3},a_{1}]=[a_{3},a_{2}]=[a_{3},a_{3}]=0.
\end{gather*}
Note also that $\mathrm{\mathbf{Leib}}(L_{6})=[L_{6},L_{6}]=\zeta^{\mathrm{right}}(L_{6})=\zeta^{\mathrm{left}}(L_{6})=\zeta(L_{6})=Fa_{3}$, $\mathrm{\mathbf{ncl}}(L_{6})=2$.

Suppose now that $[b,b]\neq0$ for every element $b\not\in\langle a_{1}\rangle$. In particular, it follows that $b\not\in\zeta(L)$. Put $[b,b]=\beta a_{3}$ where $\beta\in F$. We have $[b,a_{1}]=\gamma a_{3}$ and $[a_{1},b]=\alpha a_{3}$ for some elements $\alpha,\gamma\in F$. If $\alpha=\gamma=0$, then put $a_{2}=b$ and denote by $C$ the subalgebra generates by $a_{2}$. Then $C$ is a cyclic nilpotent ideal such that $[A,C]=[C,A]=\langle0\rangle$. Furthermore, let $u=\lambda a_{1}+\mu a_{2}+\nu a_{3}$ be the arbitrary element of $L$. Then,
\begin{gather*}
[\lambda a_{1}+\mu a_{2}+\nu a_{3},\lambda a_{1}+\mu a_{2}+\nu a_{3}]=\lambda^{2}[a_{1},a_{1}]+\mu^{2}[a_{2},a_{2}]=\\
\lambda^{2}a_{3}+\beta\mu^{2}a_{3}=(\lambda^{2}+\beta\mu^{2})a_{3}.
\end{gather*}
If $u\not\in\langle a_{1}\rangle$, then $(\lambda,\mu)\neq(0,0)$. It follows that polynomial $X^{2}+\beta$ has no root in field $F$. Thus, we come to the following type of nilpotent Leibniz algebras:
\begin{gather*}
L_{7}=Fa_{1}\oplus Fa_{2}\oplus Fa_{3}\ \mbox{where }[a_{1},a_{1}]=a_{3},[a_{2},a_{2}]=\beta a_{3}\ (\beta\neq0),\\
[a_{1},a_{2}]=[a_{1},a_{3}]=[a_{2},a_{1}]=[a_{2},a_{3}]=[a_{3},a_{1}]=[a_{3},a_{2}]=[a_{3},a_{3}]=0.
\end{gather*}
Note also that $\mathrm{\mathbf{Leib}}(L_{7})=[L_{7},L_{7}]=\zeta^{\mathrm{right}}(L_{7})=\zeta^{\mathrm{left}}(L_{7})=\zeta(L_{7})=Fa_{3}$. Moreover, polynomial $X^{2}+\beta$ has no root in field $F$, $\mathrm{\mathbf{ncl}}(L_{7})=2$.

Suppose now that $\gamma=0$ and $\alpha\neq0$. Put $a_{2}=b$. Then, $[a_{2},a_{1}]=0$, $[a_{1},a_{2}]=\alpha a_{3}$, $[a_{2},a_{2}]=\beta a_{3}$. Let $\lambda a_{1}+\mu a_{2}+\nu a_{3}$ be an arbitrary element of $L$. Then,
\begin{gather*}
[\lambda a_{1}+\mu a_{2}+\nu a_{3},\lambda a_{1}+\mu a_{2}+\nu a_{3}]=\\
\lambda^{2}[a_{1},a_{1}]+\lambda\mu[a_{1},a_{2}]+\lambda\mu[a_{2},a_{1}]+\mu^{2}[a_{2},a_{2}]=\\
\lambda^{2}a_{3}+\lambda\mu\alpha a_{3}+\mu^{2}\beta a_{3}=(\lambda^{2}+\lambda\mu\alpha+\mu^{2}\beta)a_{3}.
\end{gather*}
As above, $\lambda\neq0$, $\mu\neq0$. It follows that polynomial $X^{2}+\alpha X+\beta$ has no root in field $F$. Thus, we come to the following type of nilpotent Leibniz algebras:
\begin{gather*}
L_{8}=Fa_{1}\oplus Fa_{2}\oplus Fa_{3}\ \mbox{where }[a_{1},a_{1}]=a_{3},\\
[a_{1},a_{2}]=\alpha a_{3}\ (\alpha\neq0),[a_{2},a_{2}]=\beta a_{3}\ (\beta\neq0),\\
[a_{1},a_{3}]=[a_{2},a_{1}]=[a_{2},a_{3}]=[a_{3},a_{1}]=[a_{3},a_{2}]=[a_{3},a_{3}]=0.
\end{gather*}
Note also that $\mathrm{\mathbf{Leib}}(L_{8})=[L_{8},L_{8}]=\zeta^{\mathrm{right}}(L_{8})=\zeta^{\mathrm{left}}(L_{8})=\zeta(L_{8})=Fa_{3}$. Moreover, polynomial $X^{2}+\alpha X+\beta$ has no root in field $F$, $\mathrm{\mathbf{ncl}}(L_{8})=2$.

Suppose now that $\gamma\neq0$ and $\alpha=0$. Put $a_{2}=\gamma^{-1}b$. Then, $[a_{2},a_{1}]=a_{3}$, $[a_{1},a_{2}]=0$, $[a_{2},a_{2}]=\gamma^{-2}\beta a_{3}=\sigma a_{3}$. Let $\lambda a_{1}+\mu a_{2}+\nu a_{3}$ be an arbitrary element of $L$. Then,
\begin{gather*}
[\lambda a_{1}+\mu a_{2}+\nu a_{3},\lambda a_{1}+\mu a_{2}+\nu a_{3}]=\\
\lambda^{2}[a_{1},a_{1}]+\lambda\mu[a_{2},a_{1}]+\mu^{2}[a_{2},a_{2}]=\\
\lambda^{2}a_{3}+\lambda\mu a_{3}+\mu^{2}\sigma a_{3}=(\lambda^{2}+\lambda\mu+\mu^{2}\sigma)a_{3}.
\end{gather*}
As above, $\lambda\neq0$, $\mu\neq0$. It follows that polynomial $X^{2}+X+\sigma$ has no root in field $F$. Thus, we come to the following type of nilpotent Leibniz algebras:
\begin{gather*}
L_{9}=Fa_{1}\oplus Fa_{2}\oplus Fa_{3}\ \mbox{where }[a_{1},a_{1}]=a_{3},[a_{2},a_{1}]=a_{3},\\
[a_{2},a_{2}]=\sigma a_{3}\ (\sigma\neq0),\\
[a_{1},a_{2}]=[a_{1},a_{3}]=[a_{2},a_{3}]=[a_{3},a_{1}]=[a_{3},a_{2}]=[a_{3},a_{3}]=0.
\end{gather*}
Note also that $\mathrm{\mathbf{Leib}}(L_{9})=[L_{9},L_{9}]=\zeta^{\mathrm{right}}(L_{9})=\zeta^{\mathrm{left}}(L_{9})=\zeta(L_{9})=Fa_{3}$. Moreover, polynomial $X^{2}+X+\sigma$ has no root in field $F$, $\mathrm{\mathbf{ncl}}(L_{9})=2$.

Suppose now that $\gamma\neq0$ and $\alpha\neq0$. Put $a_{2}=\gamma^{-1}b$. Then, $[a_{2},a_{1}]=a_{3}$, $[a_{1},a_{2}]=\gamma^{-1}\alpha a_{3}=\tau a_{3}$, $[a_{2},a_{2}]=\gamma^{-2}\beta a_{3}=\sigma a_{3}$. Let $\lambda a_{1}+\mu a_{2}+\nu a_{3}$ be an arbitrary element of L. Then,
\begin{gather*}
[\lambda a_{1}+\mu a_{2}+\nu a_{3},\lambda a_{1}+\mu a_{2}+\nu a_{3}]=\\
\lambda^{2}[a_{1},a_{1}]+\lambda\mu[a_{1},a_{2}]+\lambda\mu[a_{2},a_{1}]+\mu^{2}[a_{2},a_{2}]=\\
\lambda^{2}a_{3}+\lambda\mu\tau a_{3}+\lambda\mu a_{3}+\mu^{2}\sigma a_{3}=(\lambda^{2}+\lambda\mu(\tau+1)+\mu^{2}\sigma)a_{3}.
\end{gather*}
As above, $\lambda\neq0$, $\mu\neq0$. It follows that polynomial $X^{2}+(\tau+1)X+\sigma$ has no root in field $F$. Thus, we come to the following type of nilpotent Leibniz algebras:
\begin{gather*}
L_{10}=Fa_{1}\oplus Fa_{2}\oplus Fa_{3}\ \mbox{where }[a_{1},a_{1}]=a_{3},[a_{2},a_{1}]=a_{3},\\
[a_{1},a_{2}]=\tau a_{3}\ (\tau\neq0),[a_{2},a_{2}]=\sigma a_{3}\ (\sigma\neq0),\\
[a_{1},a_{3}]=[a_{2},a_{3}]=[a_{3},a_{1}]=[a_{3},a_{2}]=[a_{3},a_{3}]=0.
\end{gather*}
Note also that $\mathrm{\mathbf{Leib}}(L_{10})=[L_{10},L_{10}]=\zeta^{\mathrm{right}}(L_{10})=\zeta^{\mathrm{left}}(L_{10})=\zeta(L_{10})=Fa_{3}$. Moreover, polynomial $X^{2}+(\tau+1)X+\sigma$ has no root in field $F$, $\mathrm{\mathbf{ncl}}(L_{10})=2$.
\qed

\begin{thm}\label{T2}
Let $L$ be a Leibniz algebra of dimension $3$ over a field $F$. Suppose that $L$ is not a Lie algebra. If the center of $L$ includes the Leibniz kernel, $\mathrm{\mathbf{dim}}_{F}(\mathrm{\mathbf{Leib}}(L))=1$ and factor-algebra $L/\mathrm{\mathbf{Leib}}(L)$ is non-abelian, then $L$ is an algebra of one of the following types.
\begin{enumerate}
\item[\upshape(i)] $\mathrm{\mathbf{Lei}}_{11}(3,F)=L_{11}$ is a direct sum of ideal $A=Fa_{1}\oplus Fa_{3}$ and subalgebra $B=Fa_{2}$. Moreover, $A$ is a nilpotent cyclic Leibniz algebra of dimension $2$, $[A,B]=[B,A]=Fa_{1}$, so that $L_{11}=Fa_{1}\oplus Fa_{2}\oplus Fa_{3}$ where $[a_{1},a_{1}]=a_{3}$, $[a_{1},a_{2}]=-a_{1}$, $[a_{2},a_{1}]=a_{1}$, $[a_{1},a_{3}]=[a_{2},a_{2}]=[a_{2},a_{3}]=[a_{3},a_{1}]=[a_{3},a_{2}]=[a_{3},a_{3}]=0$, $\mathrm{\mathbf{Leib}}(L_{11})=\zeta^{\mathrm{left}}(L_{11})=\zeta^{\mathrm{right}}(L_{11})=\zeta(L_{11})=Fa_{3}$, $[L_{11},L_{11}]=Fa_{1}\oplus Fa_{3}$, $L_{11}$ is non-nilpotent.
\item[\upshape(ii)] $\mathrm{\mathbf{Lei}}_{12}(3,F)=L_{12}$ is a direct sum of ideal $A=Fa_{1}\oplus Fa_{3}$ and subalgebra $B=Fa_{2}$. Moreover, $A$ is a nilpotent cyclic Leibniz algebra of dimension $2$, $[A,B]=[B,A]=Fa_{1}\oplus Fa_{3}$, so that $L_{12}=Fa_{1}\oplus Fa_{2}\oplus Fa_{3}$ where $[a_{1},a_{1}]=a_{3}$, $[a_{1},a_{2}]=-a_{1}-\alpha a_{3}$, $[a_{2},a_{1}]=a_{1}+\alpha a_{3}\ (\alpha\neq0)$, $[a_{1},a_{3}]=[a_{2},a_{2}]=[a_{2},a_{3}]=[a_{3},a_{1}]=[a_{3},a_{2}]=[a_{3},a_{3}]=0$, $\mathrm{\mathbf{Leib}}(L_{12})=\zeta^{\mathrm{left}}(L_{12})=\zeta^{\mathrm{right}}(L_{12})=\zeta(L_{12})=Fa_{3}$, $[L_{12},L_{12}]=Fa_{1}\oplus Fa_{3}$, $L_{12}$ is non-nilpotent.
\item[\upshape(iii)] $\mathrm{\mathbf{Lei}}_{13}(3,F)=L_{13}$ is a sum of ideal $A=Fa_{1}\oplus Fa_{3}$ and subalgebra $B=Fa_{2}\oplus Fa_{3}$. Moreover, $A$, $B$ are nilpotent cyclic Leibniz algebras of dimension $2$, $[A,B]=[B,A]=Fa_{1}$, so that $L_{13}=Fa_{1}\oplus Fa_{2}\oplus Fa_{3}$ where $[a_{1},a_{1}]=a_{3}$, $[a_{1},a_{2}]=-a_{1}$, $[a_{2},a_{1}]=a_{1}$, $[a_{2},a_{2}]=\gamma a_{3}$ $(\gamma\neq0)$, $[a_{1},a_{3}]=[a_{2},a_{3}]=[a_{3},a_{1}]=[a_{3},a_{2}]=[a_{3},a_{3}]=0$, $\mathrm{\mathbf{Leib}}(L_{13})=\zeta^{\mathrm{left}}(L_{13})=\zeta^{\mathrm{right}}(L_{13})=\zeta(L_{13})=Fa_{3}$, $[L_{13},L_{13}]=Fa_{1}\oplus Fa_{3}$. Moreover, polynomial $X^{2}+\gamma$ has no root in field $F$, $L_{13}$ is non-nilpotent.
\item[\upshape(iv)] $\mathrm{\mathbf{Lei}}_{14}(3,F)=L_{14}$ is a sum of ideal $A=Fa_{1}\oplus Fa_{3}$ and subalgebra $B=Fa_{2}\oplus Fa_{3}$. Moreover, $A,B$ are nilpotent cyclic Leibniz algebras of dimension $2$, $[A,B]=[B,A]=Fa_{1}\oplus Fa_{3}$, so that $L_{14}=Fa_{1}\oplus Fa_{2}\oplus Fa_{3}$ where $[a_{1},a_{1}]=a_{3}$, $[a_{1},a_{2}]=-a_{1}-\alpha a_{3}$, $[a_{2},a_{1}]=a_{1}+\alpha a_{3}\ (\alpha\neq0)$, $[a_{2},a_{2}]=\gamma a_{3}\ (\gamma\neq0)$, $[a_{1},a_{3}]=[a_{2},a_{3}]=[a_{3},a_{1}]=[a_{3},a_{2}]=[a_{3},a_{3}]=0$, $\mathrm{\mathbf{Leib}}(L_{14})=\zeta^{\mathrm{left}}(L_{14})=\zeta^{\mathrm{right}}(L_{14})=\zeta(L_{14})=Fa_{3}$, $[L_{14},L_{14}]=Fa_{1}\oplus Fa_{3}$. Moreover, polynomial $X^{2}+\gamma$ has no root in field $F$, $L_{14}$ is non-nilpotent.
\item[\upshape(v)] $\mathrm{\mathbf{Lei}}_{15}(3,F)=L_{15}$ is a direct sum of ideal $B=Fa_{2}$ and a subalgebra $A=Fa_{1}\oplus Fa_{3}$. Moreover, $A$ is a nilpotent cyclic Leibniz algebra of dimension $2$, $[A,B]=[B,A]=Fa_{2}$, so that $L_{15}=Fa_{1}\oplus Fa_{2}\oplus Fa_{3}$ where $[a_{1},a_{1}]=a_{3}$, $[a_{1},a_{2}]=a_{2}$, $[a_{2},a_{1}]=-a_{2}$, $[a_{1},a_{3}]=[a_{2},a_{2}]=[a_{2},a_{3}]=[a_{3},a_{1}]=[a_{3},a_{2}]=[a_{3},a_{3}]=0$, $\mathrm{\mathbf{Leib}}(L_{15})=\zeta^{\mathrm{left}}(L_{15})=\zeta^{\mathrm{right}}(L_{15})=\zeta(L_{15})=Fa_{3}$, $[L_{15},L_{15}]=Fa_{2}\oplus Fa_{3}$, $L_{15}$ is non-nilpotent.
\item[\upshape(vi)] $\mathrm{\mathbf{Lei}}_{16}(3,F)=L_{16}$ is a sum of abelian ideal $B=Fa_{2}\oplus Fa_{3}$ and a subalgebra $A=Fa_{1}\oplus Fa_{3}$. Moreover, $A$ is a nilpotent cyclic Leibniz algebra of dimension $2$, $[A,B]=[B,A]=Fa_{2}\oplus Fa_{3}$, so that $L_{16}=Fa_{1}\oplus Fa_{2}\oplus Fa_{3}$ where $[a_{1},a_{1}]=a_{3}$, $[a_{1},a_{2}]=a_{2}+\alpha a_{3}$, $[a_{2},a_{1}]=-a_{2}-\alpha a_{3}\ (\alpha\neq0)$, $[a_{1},a_{3}]=[a_{2},a_{2}]=[a_{2},a_{3}]=[a_{3},a_{1}]=[a_{3},a_{2}]=[a_{3},a_{3}]=0$, $\mathrm{\mathbf{Leib}}(L_{16})=\zeta^{\mathrm{left}}(L_{16})=\zeta^{\mathrm{right}}(L_{16})=\zeta(L_{16})=Fa_{3}$, $[L_{16},L_{16}]=Fa_{2}\oplus Fa_{3}$, $L_{16}$ is non-nilpotent.
\end{enumerate}
\end{thm}
\pf We note that the center $\zeta(L)$ has dimension at most 2. Suppose first that $\mathrm{\mathbf{dim}}_{F}(\zeta(L))=2$. Since $L$ is a not Lie algebra, there is an element $a_{1}$ such that $[a_{1},a_{1}]=a_{3}\neq0$. We note that $a_{3}\in\zeta(L)$. It follows that $[a_{1},a_{3}]=[a_{3},a_{1}]=[a_{3},a_{3}]=0$. Being an abelian algebra of dimension 2, $\zeta(L)$ has a direct decomposition $\zeta(L)=Fa_{2}\oplus Fa_{3}$ for some element $a_{2}$. We have $[a_{1},a_{2}]=[a_{2},a_{1}]=0$. But, in this case, the factor-algebra $L/\mathrm{\mathbf{Leib}}(L)$ is abelian, and we obtain a contradiction. This contradiction shows that $\zeta(L)$ has dimension 1 and hence, $\zeta(L)=\mathrm{\mathbf{Leib}}(L)$.

As noted above, $L/\mathrm{\mathbf{Leib}}(L)$ has an ideal $C/\mathrm{\mathbf{Leib}}(L)$ of dimension 1 (i.e., $C=Fc\oplus\mathrm{\mathbf{Leib}}(L)$ for some element $c$). If $[c,c]\neq0$ without loss of generality, we can put $c=a_{1}$. The ideal $\langle a_{1}\rangle=Fa_{1}\oplus Fa_{3}=A$ is nilpotent and has codimension 1. Let $b$ be an element such that $L=A\oplus Fb$. We have $[b,b]=\gamma a_{3}$ for some element $\gamma\in F$. As noted above, in this case, $[b,a_{1}]\in a_{1}+Fa_{3}$ so that $[b,a_{1}]=a_{1}+\alpha a_{3}$ for some element $\alpha\in F$. We have also $[a_{1},b]=-a_{1}+\beta a_{3}$ for some element $\beta\in F$. Using the equality
$$[b,[a_{1},b]]=[[b,a_{1}],b]+[a_{1},[b,b]]=[[b,a_{1}],b],$$
we obtain $[b,-a_{1}+\beta a_{3}]=[a_{1}+\alpha a_{3},b]$. It follows that $-[b,a_{1}]=[a_{1},b]$, so that $[a_{1},b]=-a_{1}-\alpha a_{3}$.

Suppose that $\gamma=\alpha=0$. Put $a_{2}=b$. Then, $[a_{2},a_{2}]=0$, $[a_{1},a_{2}]=-a_{1}$, $[a_{2},a_{1}]=a_{1}$. Thus, we come to the following type of Leibniz algebras:
\begin{gather*}
L_{11}=Fa_{1}\oplus Fa_{2}\oplus Fa_{3}\ \mbox{where }[a_{1},a_{1}]=a_{3},[a_{1},a_{2}]=-a_{1},[a_{2},a_{1}]=a_{1},\\
[a_{1},a_{3}]=[a_{2},a_{2}]=[a_{2},a_{3}]=[a_{3},a_{1}]=[a_{3},a_{2}]=[a_{3},a_{3}]=0.
\end{gather*}
Note also that $\mathrm{\mathbf{Leib}}(L_{11})=\zeta^{\mathrm{left}}(L_{11})=\zeta^{\mathrm{right}}(L_{11})=\zeta(L_{11})=Fa_{3}$, $[L_{11},L_{11}]=Fa_{1}\oplus Fa_{3}$, $L_{11}$ is non-nilpotent.

Suppose now that $\gamma=0$ and $\alpha\neq0$. Put again $a_{2}=b$. Then, $[a_{2},a_{2}]=0$, $[a_{1},a_{2}]=-a_{1}-\alpha a_{3}$, $[a_{2},a_{1}]=a_{1}+\alpha a_{3}$. Then, we come to the following type of Leibniz algebras:
\begin{gather*}
L_{12}=Fa_{1}\oplus Fa_{2}\oplus Fa_{3}\ \mbox{where }[a_{1},a_{1}]=a_{3},\\
[a_{1},a_{2}]=-a_{1}-\alpha a_{3},[a_{2},a_{1}]=a_{1}+\alpha a_{3}\ (\alpha\neq0),\\
[a_{1},a_{3}]=[a_{2},a_{2}]=[a_{2},a_{3}]=[a_{3},a_{1}]=[a_{3},a_{2}]=[a_{3},a_{3}]=0.
\end{gather*}
Note also that $\mathrm{\mathbf{Leib}}(L_{12})=\zeta^{\mathrm{left}}(L_{12})=\zeta^{\mathrm{right}}(L_{12})=\zeta(L_{12})=Fa_{3}$, $[L_{12},L_{12}]=Fa_{1}\oplus Fa_{3}$, $L_{12}$ is non-nilpotent.

Suppose that $\gamma\neq0$ and $\alpha=0$. Put again $a_{2}=b$. Then, $[a_{2},a_{2}]=\gamma a_{3}$, $[a_{1},a_{2}]=-a_{1}$, $[a_{2},a_{1}]=a_{1}$. Let $\lambda a_{1}+\mu a_{2}+\nu a_{3}$ be an arbitrary element of $L$. Then,
\begin{gather*}
[\lambda a_{1}+\mu a_{2}+\nu a_{3},\lambda a_{1}+\mu a_{2}+\nu a_{3}]=\\
\lambda^{2}[a_{1},a_{1}]+\lambda\mu[a_{1},a_{2}]+\lambda\mu[a_{2},a_{1}]+\mu^{2}[a_{2},a_{2}]=\\
\lambda^{2}a_{3}-\lambda\mu a_{1}+\lambda\mu a_{1}+\mu^{2}\gamma a_{3}=(\lambda^{2}+\mu^{2}\gamma)a_{3}.
\end{gather*}
As above, $\lambda\neq0$, $\mu\neq0$. It follows that polynomial $X^{2}+\gamma$ has no root in field $F$. Then, we come to the following type of Leibniz algebras:
\begin{gather*}
L_{13}=Fa_{1}\oplus Fa_{2}\oplus Fa_{3}\ \mbox{where }[a_{1},a_{1}]=a_{3},[a_{1},a_{2}]=-a_{1},[a_{2},a_{1}]=a_{1},\\
[a_{2},a_{2}]=\gamma a_{3}\ (\gamma\neq0),\\
[a_{1},a_{3}]=[a_{2},a_{3}]=[a_{3},a_{1}]=[a_{3},a_{2}]=[a_{3},a_{3}]=0.
\end{gather*}
Note also that $\mathrm{\mathbf{Leib}}(L_{13})=\zeta^{\mathrm{left}}(L_{13})=\zeta^{\mathrm{right}}(L_{13})=\zeta(L_{13})=Fa_{3}$, $[L_{13},L_{13}]=Fa_{1}\oplus Fa_{3}$. Moreover, polynomial $X^{2}+\gamma$ has no root in field $F$, $L_{13}$ is non-nilpotent.

Suppose now that $\gamma\neq0$ and $\alpha\neq0$. Put again $a_{2}=b$. Then, $[a_{2},a_{2}]=\gamma a_{3}$, $[a_{1},a_{2}]=-a_{1}-\alpha a_{3}$, $[a_{2},a_{1}]=a_{1}+\alpha a_{3}$. Let $\lambda a_{1}+\mu a_{2}+\nu a_{3}$ be an arbitrary element of $L$. Then,
\begin{gather*}
[\lambda a_{1}+\mu a_{2}+\nu a_{3},\lambda a_{1}+\mu a_{2}+\nu a_{3}]=\\
\lambda^{2}[a_{1},a_{1}]+\lambda\mu[a_{1},a_{2}]+\lambda\mu[a_{2},a_{1}]+\mu^{2}[a_{2},a_{2}]=\\
\lambda^{2}a_{3}-\lambda\mu(-a_{1}-\alpha a_{3})+\lambda\mu(a_{1}+\alpha a_{3})+\mu^{2}\gamma a_{3}=(\lambda^{2}+\mu^{2}\gamma)a_{3}.
\end{gather*}
As above, $\lambda\neq0$, $\mu\neq0$. It follows that polynomial $X^{2}+\gamma$ has no root in field $F$. Then, we come to the following type of Leibniz algebras:
\begin{gather*}
L_{14}=Fa_{1}\oplus Fa_{2}\oplus Fa_{3}\ \mbox{where }[a_{1},a_{1}]=a_{3},\\
[a_{1},a_{2}]=-a_{1}-\alpha a_{3},[a_{2},a_{1}]=a_{1}+\alpha a_{3}\ (\alpha\neq0),\\
[a_{2},a_{2}]=\gamma a_{3}\ (\gamma\neq0),\\
[a_{1},a_{3}]=[a_{2},a_{3}]=[a_{3},a_{1}]=[a_{3},a_{2}]=[a_{3},a_{3}]=0.
\end{gather*}
Note also that $\mathrm{\mathbf{Leib}}(L_{14})=\zeta^{\mathrm{left}}(L_{14})=\zeta^{\mathrm{right}}(L_{14})=\zeta(L_{14})=Fa_{3}$, $[L_{14},L_{14}]=Fa_{1}\oplus Fa_{3}$. Moreover, polynomial $X^{2}+\gamma$ has no root in field $F$, $L_{14}$ is non-nilpotent.

Suppose now that $[c,c]=0$. Put again $\mathrm{\mathbf{Leib}}(L)=Fa_{3}$. Since $\mathrm{\mathbf{Leib}}(L)=\zeta(L)$, the ideal $C=Fc\oplus\mathrm{\mathbf{Leib}}(L)$ is abelian. Suppose that there exists an element $b\not\in C$ such that $[b,b]=0$. Using the above arguments without loss of generality, we can assume that $[b,c]\in c+Fa_{3}$, $[c,b]\in -c+Fa_{3}$, so that $[b,c]=c+\alpha a_{3}$, $[c,b]=-c+\beta a_{3}$ for some elements $\alpha,\beta\in F$. Using the equality
$$[[b,c],b]=[b,[c,b]]-[c,[b,b]]=[b,[c,b]],$$
we obtain $[c+\alpha a_{3},b]=[b,-c+\beta a_{3}]$. It follows that $[c,b]=-[b,c]$, so that $[c,b]=-c-\alpha a_{3}$. Let $u=\lambda c+\mu b+\nu a_{3}$ be an arbitrary element of $L$. Then,
\begin{gather*}
[\lambda c+\mu b+\nu a_{3},\lambda c+\mu b+\nu a_{3}]=\\
\lambda^{2}[c,c]+\lambda\mu[c,b]+\lambda\mu[b,c]+\mu^{2}[b,b]=\lambda\mu[c,b]+\lambda\mu[b,c]=0.
\end{gather*}
Thus, we obtain a contradiction with the fact that $L$ is not a Lie algebra. This contradiction shows that $[b,b]\neq0$ for every element $b\not\in C$. Hence, $[b,b]=\gamma a_{3}$ where $\gamma$ is a non-zero element of field $F$. Without loss of generality, we can assume that $[b,b]=a_{3}$. Since $[b,b]\in\mathrm{\mathbf{Leib}}(L)=\zeta(L)$, we obtain that $[c,b]=-[b,c]$, so that $[c,b]=-c-\alpha a_{3}$.

If $\alpha=0$, then $[b,c]=c$, $[c,b]=-c$. Put $b=a_{1}$, $c=a_{2}$, then $Fa_{2}$ is an ideal of $L$, and we come to the following type of Leibniz algebras:
\begin{gather*}
L_{15}=Fa_{1}\oplus Fa_{2}\oplus Fa_{3}\ \mbox{where }[a_{1},a_{1}]=a_{3},[a_{1},a_{2}]=a_{2},[a_{2},a_{1}]=-a_{2},\\
[a_{1},a_{3}]=[a_{2},a_{2}]=[a_{2},a_{3}]=[a_{3},a_{1}]=[a_{3},a_{2}]=[a_{3},a_{3}]=0.
\end{gather*}
Note also that $\mathrm{\mathbf{Leib}}(L_{15})=\zeta^{\mathrm{left}}(L_{15})=\zeta^{\mathrm{right}}(L_{15})=\zeta(L_{15})=Fa_{3}$, $[L_{15},L_{15}]=Fa_{2}\oplus Fa_{3}$, $L_{15}$ is non-nilpotent.

If $\alpha\neq0$, then $[b,c]=c+\alpha a_{3}$, $[c,b]=-c-\alpha a_{3}$, and we come to the following type of Leibniz algebras:
\begin{gather*}
L_{16}=Fa_{1}\oplus Fa_{2}\oplus Fa_{3}\ \mbox{where }[a_{1},a_{1}]=a_{3},\\
[a_{1},a_{2}]=a_{2}+\alpha a_{3},[a_{2},a_{1}]=-a_{2}-\alpha a_{3}\ (\alpha\neq0),\\ [a_{1},a_{3}]=[a_{2},a_{2}]=[a_{2},a_{3}]=[a_{3},a_{1}]=[a_{3},a_{2}]=[a_{3},a_{3}]=0.
\end{gather*}
Note also that $\mathrm{\mathbf{Leib}}(L_{16})=\zeta^{\mathrm{left}}(L_{16})=\zeta^{\mathrm{right}}(L_{16})=\zeta(L_{16})=Fa_{3}$, $[L_{16},L_{16}]=Fa_{2}\oplus Fa_{3}$, $L_{16}$ is non-nilpotent.
\qed

\begin{thm}\label{T3}
Let $L$ be a Leibniz algebra of dimension $3$ over a field $F$. Suppose that $L$ is not a Lie algebra. If the center of $L$ does not include the Leibniz kernel, $\mathrm{\mathbf{dim}}_{F}(\mathrm{\mathbf{Leib}}(L))=1$ and the factor-algebra $L/\mathrm{\mathbf{Leib}}(L)$ is abelian, then $L$ is an algebra of one of the following types.
\begin{enumerate}
\item[\upshape(i)] $\mathrm{\mathbf{Lei}}_{17}(3,F)=L_{17}$ is a direct sum of two ideals $A=Fa_{1}\oplus Fa_{3}$ and $B=Fa_{2}$. Moreover, $A$ is a non-nilpotent cyclic Leibniz algebra of dimension $2$ and $B=\zeta(L)$, $[A,B]=[B,A]=\langle0\rangle$, so that $L_{17}=Fa_{1}\oplus Fa_{2}\oplus Fa_{3}$ where $[a_{1},a_{1}]=[a_{1},a_{3}]=a_{3}$, $[a_{1},a_{2}]=[a_{2},a_{1}]=[a_{2},a_{2}]=[a_{2},a_{3}]=[a_{3},a_{1}]=[a_{3},a_{2}]=[a_{3},a_{3}]=0$, $\mathrm{\mathbf{Leib}}(L_{17})=[L_{17},L_{17}]=Fa_{3}$, $\zeta^{\mathrm{left}}(L_{17})=Fa_{2}\oplus Fa_{3}$, $\zeta^{\mathrm{right}}(L_{17})=\zeta(L_{17})=Fa_{2}$, $L_{17}$ is non-nilpotent.
\item[\upshape(ii)] $\mathrm{\mathbf{Lei}}_{18}(3,F)=L_{18}$ is a direct sum of ideal $A=Fa_{1}\oplus Fa_{3}$ and a subalgebra $B=Fa_{2}$. Moreover, $A$ is a non-nilpotent cyclic Leibniz algebra of dimension $2$, $[A,B]=Fa_{3}$, $[B,A]=\langle0\rangle$, so that $L_{18}=Fa_{1}\oplus Fa_{2}\oplus Fa_{3}$ where $[a_{1},a_{1}]=[a_{1},a_{2}]=[a_{1},a_{3}]=a_{3}$, $[a_{2},a_{1}]=[a_{2},a_{2}]=[a_{2},a_{3}]=[a_{3},a_{1}]=[a_{3},a_{2}]=[a_{3},a_{3}]=0$, $\mathrm{\mathbf{Leib}}(L_{18})=[L_{18},L_{18}]=Fa_{3}$, $\zeta^{\mathrm{left}}(L_{18})=Fa_{2}\oplus Fa_{3}$, $\zeta^{\mathrm{right}}(L_{18})=\zeta(L_{18})=\langle0\rangle$, $L_{18}$ is non-nilpotent.
\item[\upshape(iii)] $\mathrm{\mathbf{Lei}}_{19}(3,F)=L_{19}$ is a direct sum of ideal $A=Fa_{1}\oplus Fa_{3}$ and a subalgebra $B=Fa_{2}$. Moreover, $A$ is a non-nilpotent cyclic Leibniz algebra of dimension $2$, $[A,B]=\langle0\rangle$, $[B,A]=Fa_{3}$, so that $L_{19}=Fa_{1}\oplus Fa_{2}\oplus Fa_{3}$ where $[a_{1},a_{1}]=[a_{1},a_{3}]=[a_{2},a_{1}]=[a_{2},a_{3}]=a_{3}$, $[a_{1},a_{2}]=[a_{2},a_{2}]=[a_{3},a_{1}]=[a_{3},a_{2}]=[a_{3},a_{3}]=0$, $\mathrm{\mathbf{Leib}}(L_{19})=[L_{19},L_{19}]=\zeta^{\mathrm{left}}(L_{19})=Fa_{3}$, $\zeta^{\mathrm{right}}(L_{19})=Fa_{2}$, $\zeta(L_{19})=\langle0\rangle$, $L_{19}$ is non-nilpotent.
\item[\upshape(iv)] $\mathrm{\mathbf{Lei}}_{20}(3,F)=L_{20}$ is a sum of two ideals $A=Fa_{1}\oplus Fa_{3}$ and $B=Fa_{2}\oplus Fa_{3}$. Furthermore, $A$ is a non-nilpotent cyclic Leibniz algebra of dimension $2$, $B$ is a nilpotent cyclic Leibniz algebra of dimension $2$, $[A,B]=Fa_{3}$, $[B,A]=\langle0\rangle$, so that $L_{20}=Fa_{1}\oplus Fa_{2}\oplus Fa_{3}$ where $[a_{1},a_{1}]=[a_{1},a_{3}]=a_{3}$, $[a_{2},a_{2}]=\sigma a_{3}\ (\sigma\neq0)$, $[a_{1},a_{2}]=[a_{2},a_{1}]=[a_{2},a_{3}]=[a_{3},a_{1}]=[a_{3},a_{2}]=[a_{3},a_{3}]=0$, $\mathrm{\mathbf{Leib}}(L_{20})=[L_{20},L_{20}]=\zeta^{\mathrm{left}}(L_{20})=Fa_{3}$, $\zeta^{\mathrm{right}}(L_{20})=\zeta(L_{20})=\langle0\rangle$. Moreover, if $\alpha a_{1}+\beta a_{2}+\gamma a_{3}\not\in A$, then $\alpha^{2}+\alpha\gamma+\beta^{2}\sigma\neq0$, $L_{20}$ is non-nilpotent.
\item[\upshape(v)] $\mathrm{\mathbf{Lei}}_{21}(3,F)=L_{21}$ is a sum of two ideals $A=Fa_{1}\oplus Fa_{3}$ and $B=Fa_{2}\oplus Fa_{3}$. Furthermore, $A$ is a non-nilpotent cyclic Leibniz algebra of dimension $2$, $B$ is a nilpotent cyclic Leibniz algebra of dimension $2$, $[A,B]=Fa_{3}$, $[B,A]=\langle0\rangle$, so that $L_{21}=Fa_{1}\oplus Fa_{2}\oplus Fa_{3}$ where $[a_{1},a_{1}]=[a_{1},a_{2}]=[a_{1},a_{3}]=a_{3}$, $[a_{2},a_{2}]=\tau a_{3}\ (\tau\neq0)$, $[a_{2},a_{1}]=[a_{2},a_{3}]=[a_{3},a_{1}]=[a_{3},a_{2}]=[a_{3},a_{3}]=0$, $\mathrm{\mathbf{Leib}}(L_{21})=[L_{21},L_{21}]=\zeta^{\mathrm{left}}(L_{21})=Fa_{3}$, $\zeta^{\mathrm{right}}(L_{21})=\zeta(L_{21})=\langle0\rangle$. Moreover, if $\alpha a_{1}+\beta a_{2}+\gamma a_{3}\not\in A$, then $\alpha^{2}+\alpha\beta+\alpha\gamma+\beta^{2}\tau\neq0$, $L_{21}$ is non-nilpotent.
\item[\upshape(vi)] $\mathrm{\mathbf{Lei}}_{22}(3,F)=L_{22}$ is a sum of two ideals $A=Fa_{1}\oplus Fa_{3}$ and $B=Fa_{2}\oplus Fa_{3}$. Furthermore, $A,B$ are non-nilpotent cyclic Leibniz algebras of dimension $2$, $[A,B]=[B,A]=Fa_{3}$, so that $L_{22}=Fa_{1}\oplus Fa_{2}\oplus Fa_{3}$ where $[a_{1},a_{1}]=[a_{1},a_{3}]=[a_{2},a_{1}]=[a_{2},a_{3}]=a_{3}$, $[a_{2},a_{2}]=\tau a_{3}\ (\tau\neq0)$, $[a_{1},a_{2}]=[a_{3},a_{1}]=[a_{3},a_{2}]=[a_{3},a_{3}]=0$, $\mathrm{\mathbf{Leib}}(L_{22})=[L_{22},L_{22}]=\zeta^{\mathrm{left}}(L_{22})=Fa_{3}$, $\zeta^{\mathrm{right}}(L_{22})=\zeta(L_{22})=\langle0\rangle$. Moreover, if $\alpha a_{1}+\beta a_{2}+\gamma a_{3}\not\in A$, then $\alpha^{2}+\alpha\gamma+\alpha\beta+\beta^{2}\tau+\beta\gamma\neq0$, $L_{22}$ is non-nilpotent.
\item[\upshape(vii)] $\mathrm{\mathbf{Lei}}_{23}(3,F)=L_{23}$ is a sum of two ideals $A=Fa_{1}\oplus Fa_{3}$ and $B=Fa_{2}\oplus Fa_{3}$. Furthermore, $A,B$ are non-nilpotent cyclic Leibniz algebras of dimension $2$, $[A,B]=[B,A]=Fa_{3}$, so that $L_{23}=Fa_{1}\oplus Fa_{2}\oplus Fa_{3}$ where $[a_{1},a_{1}]=[a_{1},a_{3}]=[a_{2},a_{1}]=[a_{2},a_{3}]=a_{3}$, $[a_{1},a_{2}]=\delta a_{3}\ (\delta\neq0)$, $[a_{2},a_{2}]=\tau a_{3}\ (\tau\neq0)$, $[a_{3},a_{1}]=[a_{3},a_{2}]=[a_{3},a_{3}]=0$, $\mathrm{\mathbf{Leib}}(L_{23})=[L_{23},L_{23}]=\zeta^{\mathrm{left}}(L_{23})=Fa_{3}$, $\zeta^{\mathrm{right}}(L_{23})=\zeta(L_{23})=\langle0\rangle$. Moreover, if $\alpha a_{1}+\beta a_{2}+\gamma a_{3}\not\in A$, then $\alpha^{2}+\alpha\beta\delta+\alpha\gamma+\alpha\beta+\beta^{2}\tau+\beta\gamma\neq0$, $L_{23}$ is non-nilpotent.
\end{enumerate}
\end{thm}
\pf Since $\mathrm{\mathbf{dim}}_{F}(\mathrm{\mathbf{Leib}}(L))=1$, $\mathrm{\mathbf{Leib}}(L)\cap\zeta(L)=\langle0\rangle$. If we suppose that the center $\zeta(L)$ has dimension 2, then $L=\mathrm{\mathbf{Leib}}(L)\oplus\zeta(L)$. But, in this case, $L$ is abelian and we obtain a contradiction. Suppose now that $\mathrm{\mathbf{dim}}_{F}(\zeta(L))=1$.

Since $L$ is not a Lie algebra, there is an element $a_{1}$ such that $[a_{1},a_{1}]=a_{3}\neq0$. Then, $\mathrm{\mathbf{Leib}}(L)=Fa_{3}$. Then, $A=\langle a_{1}\rangle=Fa_{1}\oplus Fa_{3}$ is a subalgebra of $L$. It is obvious that $A\cap\zeta(L)=\langle0\rangle$, so that $L=A\oplus\zeta(L)$, $A$ is an ideal of $L$ and $B=\zeta(L)=Fa_{2}$. If we suppose that $A$ is nilpotent (that is, $[a_{1},a_{3}]=0$), then $\mathrm{\mathbf{Leib}}(L)=\mathrm{\mathbf{Leib}}(A)=\zeta(A)\leqslant\zeta(L)$, and we obtain a contradiction. Thus, $A$ is not nilpotent. As we have seen above, we can choose an element $a_{1}$ such that $[a_{1},a_{3}]=a_{3}$. Thus, we come to the following type of Leibniz algebras:
\begin{gather*}
L_{17}=Fa_{1}\oplus Fa_{2}\oplus Fa_{3}\ \mbox{where }[a_{1},a_{1}]=[a_{1},a_{3}]=a_{3},\\
[a_{1},a_{2}]=[a_{2},a_{1}]=[a_{2},a_{2}]=[a_{2},a_{3}]=[a_{3},a_{1}]=[a_{3},a_{2}]=[a_{3},a_{3}]=0.
\end{gather*}
Note also that $\mathrm{\mathbf{Leib}}(L_{17})=[L_{17},L_{17}]=Fa_{3}$, $\zeta^{\mathrm{left}}(L_{17})=Fa_{2}\oplus Fa_{3}$, $\zeta^{\mathrm{right}}(L_{17})=\zeta(L_{17})=Fa_{2}$, $L_{17}$ is non-nilpotent.

Suppose now that $\zeta(L)$ is zero. Since $L$ is not a Lie algebra, there is an element $a_{1}$ such that $[a_{1},a_{1}]=a_{3}\neq0$. Then, $\mathrm{\mathbf{Leib}}(L)=Fa_{3}$. Since $L/\mathrm{\mathbf{Leib}}(L)$ is abelian, a subalgebra $A=\langle a_{1}\rangle=Fa_{1}\oplus Fa_{3}$ is an ideal of $L$. Suppose that $A$ is nilpotent (that is, $[a_{1},a_{3}]=0$). Let $b$ be an element of $L$ such that $b\not\in A$. We have:
$$[b,a_{3}]=[b,[a_{1},a_{1}]]=[[b,a_{1}],a_{1}]+[a_{1},[b,a_{1}]]=[\lambda a_{3},a_{1}]+[a_{1},\lambda a_{3}]=0.$$
Since $[a_{3},b]=0$, we obtain that $a_{3}\in\zeta(L)$, and we obtain a contradiction. This contradiction shows that subalgebra $A=\langle a_{1}\rangle=Fa_{1}\oplus Fa_{3}$ is not nilpotent. As we have seen above, we can choose an element $a_{1}$ such that $[a_{1},a_{3}]=a_{3}$.

Suppose that $L$ contains an element $b$ such that $b\not\in A$ and $[b,b]=0$. Since $L/\mathrm{\mathbf{Leib}}(L)$ is abelian, $[b,a_{1}]=\lambda a_{3}$, $[a_{1},b]=\mu a_{3}$ for some elements $\lambda,\mu\in F$. As above,
$$[b,a_{3}]=[a_{1},[b,a_{1}]]=[a_{1},\lambda a_{3}]=\lambda[a_{1},a_{3}]=\lambda a_{3}.$$

If $\lambda=\mu=0$, then $b\in\zeta(L)$, and we obtain a contradiction with our assumption concerning the inclusion of $\zeta(L)$.

Suppose, now, that $\lambda=0$, $\mu\neq0$. Put $a_{2}=\mu^{-1}b$, then $[a_{2},a_{1}]=0$, $[a_{1},a_{2}]=a_{3}$, $[a_{2},a_{2}]=0$. As seen above, we can observe that $[a_{2},a_{3}]=0$, and we obtain that $[a_{2},A]=\langle0\rangle$. It follows that $a_{2}\in\zeta^{\mathrm{left}}(L)$. Thus, we come to the following type of Leibniz algebras:
\begin{gather*}
L_{18}=Fa_{1}\oplus Fa_{2}\oplus Fa_{3}\ \mbox{where }[a_{1},a_{1}]=[a_{1},a_{2}]=[a_{1},a_{3}]=a_{3},\\
[a_{2},a_{1}]=[a_{2},a_{2}]=[a_{2},a_{3}]=[a_{3},a_{1}]=[a_{3},a_{2}]=[a_{3},a_{3}]=0.
\end{gather*}
Note also that $\mathrm{\mathbf{Leib}}(L_{18})=[L_{18},L_{18}]=Fa_{3}$, $\zeta^{\mathrm{left}}(L_{18})=Fa_{2}\oplus Fa_{3}$, $\zeta^{\mathrm{right}}(L_{18})=\zeta(L_{14})=\langle0\rangle$, $L_{18}$ is non-nilpotent.

Suppose now that $\lambda\neq0$, $\mu=0$. Put $a_{2}=\lambda^{-1}b$, then $[a_{2},a_{1}]=a_{3}$, $[a_{1},a_{2}]=0$, $[a_{2},a_{2}]=0$. Since $[a_{3},a_{2}]=0$, $[A,a_{2}]=\langle0\rangle$. It follows that $a_{2}\in\zeta^{\mathrm{right}}(L)$. Thus, we come to the following type of Leibniz algebras:
\begin{gather*}
L_{19}=Fa_{1}\oplus Fa_{2}\oplus Fa_{3}\ \mbox{where}\\
[a_{1},a_{1}]=[a_{1},a_{3}]=[a_{2},a_{1}]=[a_{2},a_{3}]=a_{3},\\
[a_{1},a_{2}]=[a_{2},a_{2}]=[a_{3},a_{1}]=[a_{3},a_{2}]=[a_{3},a_{3}]=0.
\end{gather*}
Note also that $\mathrm{\mathbf{Leib}}(L_{19})=[L_{19},L_{19}]=\zeta^{\mathrm{left}}(L_{19})=Fa_{3}$, $\zeta^{\mathrm{right}}(L_{19})=Fa_{2}$, $\zeta(L_{19})=\langle0\rangle$, $L_{19}$ is non-nilpotent.

Suppose now that $\mu\neq0$, $\lambda\neq0$. We have:
\begin{gather*}
0=[a_{1},0]=[a_{1},[b,b]]=[[a_{1},b],b]+[b,[a_{1},b]]=\\
[b,[a_{1},b]]=[b,\mu a_{3}]=\mu[b,a_{3}]=\mu\lambda a_{3}.
\end{gather*}
It follows that $\mu\lambda=0$, and we obtain a contradiction.

Suppose that $[b,b]\neq0$ for every element $b$ such that $b\not\in A$. Since $L/\mathrm{\mathbf{Leib}}(L)$ is abelian, $[b,a_{1}]=\lambda a_{3}$, $[a_{1},b]=\mu a_{3}$ for some elements $\lambda,\mu\in F$, and $[b,b]=\sigma a_{3}$ for some non-zero element $\sigma\in F$.

As above,
$$[b,a_{3}]=[a_{1},[b,a_{1}]]=[a_{1},\lambda a_{3}]=\lambda[a_{1},a_{3}]=\lambda a_{3}.$$

If $\lambda=\mu=0$, then $[b,a_{1}]=[a_{1},b]=[b,a_{3}]=0$. Put $a_{2}=b$. Let $u=\alpha a_{1}+\beta a_{2}+\gamma a_{3}$ be an arbitrary element of $L$. Then,
\begin{gather*}
[u,u]=[\alpha a_{1}+\beta a_{2}+\gamma a_{3},\alpha a_{1}+\beta a_{2}+\gamma a_{3}]=\\
\alpha^{2}[a_{1},a_{1}]+\alpha\gamma[a_{1},a_{3}]+\beta^{2}[a_{2},a_{2}]=\\
\alpha^{2}a_{3}+\alpha\gamma a_{3}+\beta^{2}a_{3}=(\alpha^{2}+\alpha\gamma+\beta^{2}\sigma)a_{3}.
\end{gather*}
If $u\not\in A$, then $[u,u]\neq0$. It follows that $\alpha^{2}+\alpha\gamma+\beta^{2}\sigma\neq0$. Thus, we come to the following type of Leibniz algebras:
\begin{gather*}
L_{20}=Fa_{1}\oplus Fa_{2}\oplus Fa_{3}\ \mbox{where }[a_{1},a_{1}]=[a_{1},a_{3}]=a_{3},\\
[a_{2},a_{2}]=\sigma a_{3}\ (\sigma\neq0),\\
[a_{1},a_{2}]=[a_{2},a_{1}]=[a_{2},a_{3}]=[a_{3},a_{1}]=[a_{3},a_{2}]=[a_{3},a_{3}]=0.
\end{gather*}
Note also that $\mathrm{\mathbf{Leib}}(L_{20})=[L_{20},L_{20}]=\zeta^{\mathrm{left}}(L_{20})=Fa_{3}$, $\zeta^{\mathrm{right}}(L_{20})=\zeta(L_{20})=\langle0\rangle$. Moreover, if $\alpha a_{1}+\beta a_{2}+\gamma a_{3}\not\in A$, then $\alpha^{2}+\alpha\gamma+\beta^{2}\sigma\neq0$, $L_{20}$ is non-nilpotent.

Suppose now that $\lambda=0$, $\mu\neq0$. Put $a_{2}=\mu^{-1}b$. Then $[a_{2},a_{1}]=[a_{2},a_{3}]=0$, $[a_{1},a_{2}]=a_{3}$, $[a_{2},a_{2}]=\mu^{-2}[b,b]=\mu^{-2}\sigma a_{3}=\tau a_{3}$. Let $\alpha a_{1}+\beta a_{2}+\gamma a_{3}$ be an arbitrary element of $L$. Then,
\begin{gather*}
[\alpha a_{1}+\beta a_{2}+\gamma a_{3},\alpha a_{1}+\beta a_{2}+\gamma a_{3}]=\\
\alpha^{2}[a_{1},a_{1}]+\alpha\beta[a_{1},a_{2}]+\alpha\gamma[a_{1},a_{3}]+\beta^{2}[a_{2},a_{2}]=\\
\alpha^{2}a_{3}+\alpha\beta a_{3}+\alpha\gamma a_{3}+\beta^{2}\tau a_{3}=(\alpha^{2}+\alpha\beta+\alpha\gamma+\beta^{2}\tau)a_{3}.
\end{gather*}
As above, $\alpha^{2}+\alpha\beta+\alpha\gamma+\beta^{2}\tau\neq0$. Thus, we come to the following type of Leibniz algebras:
\begin{gather*}
L_{21}=Fa_{1}\oplus Fa_{2}\oplus Fa_{3}\ \mbox{where }[a_{1},a_{1}]=[a_{1},a_{2}]=[a_{1},a_{3}]=a_{3},\\
[a_{2},a_{2}]=\tau a_{3}\ (\tau\neq0),\\
[a_{2},a_{1}]=[a_{2},a_{3}]=[a_{3},a_{1}]=[a_{3},a_{2}]=[a_{3},a_{3}]=0.
\end{gather*}
Note also that $\mathrm{\mathbf{Leib}}(L_{21})=[L_{21},L_{21}]=\zeta^{\mathrm{left}}(L_{21})=Fa_{3}$, $\zeta^{\mathrm{right}}(L_{21})=\zeta(L_{21})=\langle0\rangle$. Moreover, if $\alpha a_{1}+\beta a_{2}+\gamma a_{3}\not\in A$, then $\alpha^{2}+\alpha\beta+\alpha\gamma+\beta^{2}\tau\neq0$, $L_{21}$ is non-nilpotent.

Suppose now that $\lambda\neq0$, $\mu=0$. Put $a_{2}=\lambda^{-1}b$. Then $[a_{2},a_{1}]=[a_{2},a_{3}]=a_{3}$, $[a_{1},a_{2}]=0$, $[a_{2},a_{2}]=\lambda^{-2}[b,b]=\lambda^{-2}\sigma a_{3}=\tau a_{3}$. Let $\alpha a_{1}+\beta a_{2}+\gamma a_{3}$ be an arbitrary element of $L$. Then,
\begin{gather*}
[\alpha a_{1}+\beta a_{2}+\gamma a_{3},\alpha a_{1}+\beta a_{2}+\gamma a_{3}]=\\
\alpha^{2}[a_{1},a_{1}]+\alpha\gamma[a_{1},a_{3}]+\alpha\beta[a_{2},a_{1}]+\beta^{2}[a_{2},a_{2}]+\beta\gamma[a_{2},a_{3}]=\\
\alpha^{2}a_{3}+\alpha\gamma a_{3}+\alpha\beta a_{3}+\beta^{2}\tau a_{3}+\beta\gamma a_{3}=\\
(\alpha^{2}+\alpha\gamma+\alpha\beta+\beta^{2}\tau+\beta\gamma)a_{3}.
\end{gather*}
As above, $\alpha^{2}+\alpha\gamma+\alpha\beta+\beta^{2}\tau+\beta\gamma\neq0$. Thus, we come to the following type of Leibniz algebras:
\begin{gather*}
L_{22}=Fa_{1}\oplus Fa_{2}\oplus Fa_{3}\ \mbox{where}\\
[a_{1},a_{1}]=[a_{1},a_{3}]=[a_{2},a_{1}]=[a_{2},a_{3}]=a_{3},\\
[a_{2},a_{2}]=\tau a_{3}\ (\tau\neq0),\\
[a_{1},a_{2}]=[a_{3},a_{1}]=[a_{3},a_{2}]=[a_{3},a_{3}]=0.
\end{gather*}
Note also that $\mathrm{\mathbf{Leib}}(L_{22})=[L_{22},L_{22}]=\zeta^{\mathrm{left}}(L_{22})=Fa_{3}$, $\zeta^{\mathrm{right}}(L_{22})=\zeta(L_{22})=\langle0\rangle$. Moreover, if $\alpha a_{1}+\beta a_{2}+\gamma a_{3}\not in A$, then $\alpha^{2}+\alpha\gamma+\alpha\beta+\beta^{2}\tau+\beta\gamma\neq0$, $L_{22}$ is non-nilpotent.

Suppose now that $\lambda\neq0$, $\mu\neq0$. Put $a_{2}=\lambda^{-1}b$. Then $[a_{2},a_{1}]=[a_{2},a_{3}]=a_{3}$, $[a_{1},a_{2}]=\lambda^{-1}[a_{1},b]=\lambda^{-1}\mu a_{3}=\delta a_{3}$, $[a_{2},a_{2}]=\lambda^{-2}[b,b]=\lambda^{-2}\sigma a_{3}=\tau a_{3}$. Let $\alpha a_{1}+\beta a_{2}+\gamma a_{3}$ be an arbitrary element of $L$. Then,
\begin{gather*}
[\alpha a_{1}+\beta a_{2}+\gamma a_{3},\alpha a_{1}+\beta a_{2}+\gamma a_{3}]=\\
\alpha^{2}[a_{1},a_{1}]+\alpha\beta[a_{1},a_{2}]+\alpha\gamma[a_{1},a_{3}]+\\
\alpha\beta[a_{2},a_{1}]+\beta^{2}[a_{2},a_{2}]+\beta\gamma[a_{2},a_{3}]=\\
\alpha^{2}a_{3}+\alpha\beta\delta a_{3}+\alpha\gamma a_{3}+\alpha\beta a_{3}+\beta^{2}\tau a_{3}+\beta\gamma a_{3}=\\
(\alpha^{2}+\alpha\beta\delta+\alpha\gamma+\alpha\beta+\beta^{2}\tau+\beta\gamma)a_{3}.
\end{gather*}
As above, $\alpha^{2}+\alpha\beta\delta+\alpha\gamma+\alpha\beta+\beta^{2}\tau+\beta\gamma\neq0$. Thus, we come to the following type of Leibniz algebras:
\begin{gather*}
L_{23}=Fa_{1}\oplus Fa_{2}\oplus Fa_{3}\ \mbox{where}\\
[a_{1},a_{1}]=[a_{1},a_{3}]=[a_{2},a_{1}]=[a_{2},a_{3}]=a_{3},\\
[a_{1},a_{2}]=\delta a_{3}\ (\delta\neq0),[a_{2},a_{2}]=\tau a_{3}\ (\tau\neq0),\\
[a_{3},a_{1}]=[a_{3},a_{2}]=[a_{3},a_{3}]=0.
\end{gather*}
Note also that $\mathrm{\mathbf{Leib}}(L_{23})=[L_{23},L_{23}]=\zeta^{\mathrm{left}}(L_{23})=Fa_{3}$, $\zeta^{\mathrm{right}}(L_{23})=\zeta(L_{23})=\langle0\rangle$. Moreover, if $\alpha a_{1}+\beta a_{2}+\gamma a_{3}\not\in A$, then $\alpha^{2}+\alpha\beta\delta+\alpha\gamma+\alpha\beta+\beta^{2}\tau+\beta\gamma\neq0$, $L_{23}$ is non-nilpotent.
\qed

Let $L$ be a Leibniz algebra over a field $F$, $M$ be non-empty subset of $L$ and $H$ be a subalgebra of $L$. Put
\begin{gather*}
\mathrm{\mathbf{Ann}}_{H}^{\mathrm{left}}(M)=\{a\in H|\ [a,M]=\langle0\rangle\},\\
\mathrm{\mathbf{Ann}}_{H}^{\mathrm{right}}(M)=\{a\in H|\ [M,a]=\langle0\rangle\}.
\end{gather*}
The subset $\mathrm{\mathbf{Ann}}_{H}^{\mathrm{left}}(M)$ is called the \textit{left annihilator} of $M$ in subalgebra $H$. The subset $\mathrm{\mathbf{Ann}}_{H}^{\mathrm{right}}(M)$ is called the \textit{right annihilator} of $M$ in subalgebra $H$. The intersection
$$\mathrm{\mathbf{Ann}}_{H}(M)=\mathrm{\mathbf{Ann}}_{H}^{\mathrm{left}}(M)\cap\mathrm{\mathbf{Ann}}_{H}^{\mathrm{right}}(M)$$
is called the \textit{annihilator} of $M$ in subalgebra $H$.

It is not hard to see that all of these subsets are subalgebras of $L$. Moreover, if $M$ is an ideal of $L$, then $\mathrm{\mathbf{Ann}}_{H}(M)$ is an ideal of $L$ (see, for example,~\cite{KOP2016}).

\begin{thm}\label{T4}
Let $L$ be a Leibniz algebra over a field $F$ having dimension $3$. Suppose that $L$ is not a Lie algebra. Suppose that the center of $L$ does not include the Leibniz kernel, $\mathrm{\mathbf{dim}}_{F}(\mathrm{\mathbf{Leib}}(L))=1$, and the factor-algebra $L/\mathrm{\mathbf{Leib}}(L)$ is non-abelian. Then, $L$ is an algebra of one of the following types.
\begin{enumerate}
\item[\upshape(i)] $\mathrm{\mathbf{Lei}}_{24}(3,F)=L_{24}$ is a direct sum of ideal $A=Fa_{1}\oplus Fa_{3}$ and subalgebra $B=Fa_{2}$. Moreover, $A$ is a nilpotent cyclic Leibniz algebra of dimension $2$, $\mathrm{\mathbf{char}}(F)\neq2$, $[A,B]=Fa_{1}$, $[B,A]=Fa_{1}\oplus Fa_{3}$, so that $L_{24}=Fa_{1}\oplus Fa_{2}\oplus Fa_{3}$ where $[a_{1},a_{1}]=a_{3}$, $[a_{1},a_{2}]=-a_{1}$, $[a_{2},a_{1}]=a_{1}$, $[a_{2},a_{3}]=2a_{3}$, $[a_{1},a_{3}]=[a_{2},a_{2}]=[a_{3},a_{1}]=[a_{3},a_{2}]=[a_{3},a_{3}]=0$, $\mathrm{\mathbf{Leib}}(L_{24})=\zeta^{\mathrm{left}}(L_{24})=Fa_{3}$, $\zeta^{\mathrm{right}}(L_{24})=\zeta(L_{24})=\langle0\rangle$, $[L_{24},L_{24}]=Fa_{1}\oplus Fa_{3}$, $L_{24}$ is non-nilpotent.
\item[\upshape(ii)] $\mathrm{\mathbf{Lei}}_{25}(3,F)=L_{25}$ is a direct sum of ideal $A=Fa_{1}\oplus Fa_{3}$ and subalgebra $B=Fa_{2}$. Moreover, $A$ is a nilpotent cyclic Leibniz algebra of dimension $2$, $\mathrm{\mathbf{char}}(F)\neq2$, $[A,B]=[B,A]=Fa_{1}\oplus Fa_{3}$, so that $L_{25}=Fa_{1}\oplus Fa_{2}\oplus Fa_{3}$ where $[a_{1},a_{1}]=a_{3}$, $[a_{1},a_{2}]=-a_{1}+\alpha a_{3}$, $[a_{2},a_{1}]=a_{1}+\alpha a_{3}\ (\alpha\neq0)$, $[a_{2},a_{3}]=2a_{3}$, $[a_{1},a_{3}]=[a_{2},a_{2}]=[a_{3},a_{1}]=[a_{3},a_{2}]=[a_{3},a_{3}]=0$, $\mathrm{\mathbf{Leib}}(L_{25})=\zeta^{\mathrm{left}}(L_{25})=Fa_{3}$, $\zeta^{\mathrm{right}}(L_{25})=\zeta(L_{25})=\langle0\rangle$, $[L_{25},L_{25}]=Fa_{1}\oplus Fa_{3}$, $L_{25}$ is non-nilpotent.
\item[\upshape(iii)] $\mathrm{\mathbf{Lei}}_{26}(3,F)=L_{26}$ is a direct sum of ideal $A=Fa_{1}\oplus Fa_{3}$ and subalgebra $B=Fa_{2}$. Furthermore, $A$ is a non-nilpotent cyclic Leibniz algebra of dimension $2$, $[A,B]=Fa_{1}$, $[B,A]=Fa_{1}\oplus Fa_{3}$, so that $L_{26}=Fa_{1}\oplus Fa_{2}\oplus Fa_{3}$ where $[a_{1},a_{1}]=[a_{1},a_{3}]=a_{3}$, $[a_{1},a_{2}]=-a_{1}$, $[a_{2},a_{1}]=a_{1}$, $[a_{2},a_{3}]=2a_{3}$, $[a_{2},a_{2}]=[a_{3},a_{1}]=[a_{3},a_{2}]=[a_{3},a_{3}]=0$, $\mathrm{\mathbf{Leib}}(L_{26})=\zeta^{\mathrm{left}}(L_{26})=Fa_{3}$, $[L_{26},L_{26}]=Fa_{1}\oplus Fa_{3}$, $\zeta^{\mathrm{right}}(L_{26})=\zeta(L_{26})=\langle0\rangle$. Moreover, $\mathrm{\mathbf{char}}(F)\neq2$, $L_{26}$ is non-nilpotent.
\item[\upshape(iv)] $\mathrm{\mathbf{Lei}}_{27}(3,F)=L_{27}$ is a direct sum of ideal $A=Fa_{1}\oplus Fa_{3}$ and subalgebra $B=Fa_{2}$. Moreover, $A$ is a non-nilpotent cyclic Leibniz algebra of dimension $2$, $[A,B]=[B,A]=Fa_{1}\oplus Fa_{3}$, so that $L_{27}=Fa_{1}\oplus Fa_{2}\oplus Fa_{3}$ where $[a_{1},a_{1}]=[a_{1},a_{3}]=a_{3}$, $[a_{1},a_{2}]=-a_{1}+\beta a_{3}$ $(\beta=\alpha(1+\alpha)^{-1})$, $[a_{2},a_{1}]=a_{1}+\alpha a_{3}$, $[a_{2},a_{3}]=(2+\alpha)a_{3}$ $(\alpha\neq0,\alpha\neq-1,\alpha\neq-2)$, $[a_{2},a_{2}]=[a_{3},a_{1}]=[a_{3},a_{2}]=[a_{3},a_{3}]=0$, $\mathrm{\mathbf{Leib}}(L_{27})=\zeta^{\mathrm{left}}(L_{27})=Fa_{3}$, $[L_{27},L_{27}]=Fa_{1}\oplus Fa_{3}$, $\zeta^{\mathrm{right}}(L_{27})=\zeta(L_{27})=\langle0\rangle$, $L_{27}$ is non-nilpotent.
\item[\upshape(v)] $\mathrm{\mathbf{Lei}}_{28}(3,F)=L_{28}$ is a sum of ideal $A=Fa_{1}\oplus Fa_{3}$ and subalgebra $B=Fa_{2}\oplus Fa_{3}$. Furthermore, $A$ is a nilpotent cyclic Leibniz algebra of dimension $2$, $B$ is a non-nilpotent cyclic Leibniz algebra of dimension $2$, $\mathrm{\mathbf{char}}(F)\neq2$, $[A,B]=Fa_{1}$, $[B,A]=Fa_{1}\oplus Fa_{3}$, so that $L_{28}=Fa_{1}\oplus Fa_{2}\oplus Fa_{3}$ where $[a_{1},a_{1}]=a_{3}$, $[a_{1},a_{2}]=-a_{1}$, $[a_{2},a_{1}]=a_{1}$, $[a_{2},a_{2}]=\gamma a_{3}$ $(\gamma\neq0)$, $[a_{2},a_{3}]=2a_{3}$, $[a_{1},a_{3}]=[a_{3},a_{1}]=[a_{3},a_{2}]=[a_{3},a_{3}]=0$, $\mathrm{\mathbf{Leib}}(L_{28})=\zeta^{\mathrm{left}}(L_{28})=Fa_{3}$, $[L_{28},L_{28}]=Fa_{1}\oplus Fa_{3}$, $\zeta^{\mathrm{right}}(L_{28})=\zeta(L_{28})=\langle0\rangle$. Moreover, if $\lambda a_{1}+\mu a_{2}+\nu a_{3}\not\in A$, then $\lambda^{2}+\mu^{2}\gamma+2\mu\nu\neq0$, $L_{28}$ is non-nilpotent.
\item[\upshape(vi)] $\mathrm{\mathbf{Lei}}_{29}(3,F)=L_{29}$ is a sum of ideal $A=Fa_{1}\oplus Fa_{3}$ and subalgebra $B=Fa_{2}\oplus Fa_{3}$. Furthermore, $A$ is a nilpotent cyclic Leibniz algebra of dimension $2$, $B$ is a non-nilpotent cyclic Leibniz algebra of dimension $2$, $\mathrm{\mathbf{char}}(F)\neq2$, $[A,B]=[B,A]=Fa_{1}\oplus Fa_{3}$, so that $L_{29}=Fa_{1}\oplus Fa_{2}\oplus Fa_{3}$ where $[a_{1},a_{1}]=a_{3}$, $[a_{1},a_{2}]=-a_{1}+\alpha a_{3}$, $[a_{2},a_{1}]=a_{1}+\alpha a_{3}\ (\alpha\neq0)$, $[a_{2},a_{2}]=\gamma a_{3}\ (\gamma\neq0)$, $[a_{2},a_{3}]=2a_{3}$, $[a_{1},a_{3}]=[a_{3},a_{1}]=[a_{3},a_{2}]=[a_{3},a_{3}]=0$, $\mathrm{\mathbf{Leib}}(L_{29})=\zeta^{\mathrm{left}}(L_{29})=Fa_{3}$, $[L_{29},L_{29}]=Fa_{1}\oplus Fa_{3}$, $\zeta^{\mathrm{right}}(L_{29})=\zeta(L_{29})=\langle0\rangle$. Moreover, if $\lambda a_{1}+\mu a_{2}+\nu a_{3}\not\in A$, then $\lambda^{2}+2\lambda\mu\alpha+\mu^{2}\gamma+2\mu\nu\neq0$, $L_{29}$ is non-nilpotent.
\item[\upshape(vii)] $\mathrm{\mathbf{Lei}}_{30}(3,F)=L_{30}$ is a sum of ideal $A=Fa_{1}\oplus Fa_{3}$ and subalgebra $B=Fa_{2}\oplus Fa_{3}$. Moreover, $A$, $B$ are non-nilpotent cyclic Leibniz algebras of dimension $2$, $\mathrm{\mathbf{char}}(F)\neq2$, $[A,B]=[B,A]=Fa_{1}\oplus Fa_{3}$, so that $L_{30}=Fa_{1}\oplus Fa_{2}\oplus Fa_{3}$ where $[a_{1},a_{1}]=[a_{1},a_{3}]=a_{3}$, $[a_{1},a_{2}]=-a_{1}+\gamma a_{3}$ $(\gamma\neq0)$, $[a_{2},a_{1}]=a_{1}$, $[a_{2},a_{2}]=\gamma a_{3}$, $[a_{2},a_{3}]=2a_{3}$, $[a_{3},a_{1}]=[a_{3},a_{2}]=[a_{3},a_{3}]=0$, $\mathrm{\mathbf{Leib}}(L_{30})=\zeta^{\mathrm{left}}(L_{30})=Fa_{3}$, $[L_{30},L_{30}]=Fa_{1}\oplus Fa_{3}$, $\zeta^{\mathrm{right}}(L_{30})=\zeta(L_{30})=\langle0\rangle$. Moreover, if $\lambda a_{1}+\mu a_{2}+\nu a_{3}\not\in A$, then $\lambda^{2}+\lambda\mu\gamma+\lambda\nu+\mu^{2}\gamma+2\mu\nu\neq0$, $L_{30}$ is non-nilpotent.
\item[\upshape(viii)] $\mathrm{\mathbf{Lei}}_{31}(3,F)=L_{31}$ is a sum of ideal $A=Fa_{1}\oplus Fa_{3}$ and subalgebra $B=Fa_{2}\oplus Fa_{3}$. Moreover, $A,B$ are non-nilpotent cyclic Leibniz algebras of dimension $2$, $[A,B]=[B,A]=Fa_{1}\oplus Fa_{3}$, so that $L_{31}=Fa_{1}\oplus Fa_{2}\oplus Fa_{3}$ where $[a_{1},a_{1}]=[a_{1},a_{3}]=[a_{2},a_{2}]=[a_{2},a_{3}]=a_{3}$, $[a_{1},a_{2}]=-a_{1}$, $[a_{2},a_{1}]=a_{1}-a_{3}$, $[a_{3},a_{1}]=[a_{3},a_{2}]=[a_{3},a_{3}]=0$, $\mathrm{\mathbf{Leib}}(L_{31})=\zeta^{\mathrm{left}}(L_{31})=Fa_{3}$, $[L_{31},L_{31}]=Fa_{1}\oplus Fa_{3}$, $\zeta^{\mathrm{right}}(L_{31})=\zeta(L_{31})=\langle0\rangle$. Moreover, if $\lambda a_{1}+\mu a_{2}+\nu a_{3}\not\in A$, then $\lambda^{2}+\lambda\nu-\lambda\mu+\mu^{2}+\mu\nu\neq0$, $L_{31}$ is non-nilpotent.
\item[\upshape(ix)] $\mathrm{\mathbf{Lei}}_{32}(3,F)=L_{31}$ is a sum of ideal $A=Fa_{1}\oplus Fa_{3}$ and subalgebra $B=Fa_{2}\oplus Fa_{3}$. Moreover, $A,B$ are non-nilpotent cyclic Leibniz algebras of dimension $2$, $[A,B]=[B,A]=Fa_{1}\oplus Fa_{3}$, so that $L_{32}=Fa_{1}\oplus Fa_{2}\oplus Fa_{3}$ where $[a_{1},a_{1}]=[a_{1},a_{3}]=[a_{2},a_{2}]=[a_{2},a_{3}]=a_{3}$, $[a_{1},a_{2}]=-a_{1}+\beta a_{3}\ (\beta\neq0)$, $[a_{2},a_{1}]=a_{1}-a_{3}$, $[a_{3},a_{1}]=[a_{3},a_{2}]=[a_{3},a_{3}]=0$, $\mathrm{\mathbf{Leib}}(L_{32})=\zeta^{\mathrm{left}}(L_{32})=Fa_{3}$, $[L_{32},L_{32}]=Fa_{1}\oplus Fa_{3}$, $\zeta^{\mathrm{right}}(L_{32})=\zeta(L_{32})=\langle0\rangle$. Moreover, if $\lambda a_{1}+\mu a_{2}+\nu a_{3}\not\in A$, then $\lambda^{2}+\lambda\mu\beta+\lambda\nu-\lambda\mu+\mu^{2}+\mu\nu\neq0$, $L_{32}$ is non-nilpotent.
\item[\upshape(x)] $\mathrm{\mathbf{Lei}}_{33}(3,F)=L_{31}$ is a sum of ideal $A=Fa_{1}\oplus Fa_{3}$ and subalgebra $B=Fa_{2}\oplus Fa_{3}$. Moreover, $A,B$ are non-nilpotent cyclic Leibniz algebras of dimension $2$, $[A,B]=[B,A]=Fa_{1}\oplus Fa_{3}$, $L_{33}=Fa_{1}\oplus Fa_{2}\oplus Fa_{3}$ where $[a_{1},a_{1}]=[a_{1},a_{3}]=a_{3}$, $[a_{1},a_{2}]=-a_{1}$, $[a_{2},a_{1}]=a_{1}-\gamma a_{3}\ (\gamma\neq0,\gamma\neq1,\gamma\neq2)$, $[a_{2},a_{2}]=\gamma a_{3}$, $[a_{2},a_{3}]=(2-\gamma)a_{3}$, $[a_{3},a_{1}]=[a_{3},a_{2}]=[a_{3},a_{3}]=0$, $\mathrm{\mathbf{Leib}}(L_{33})=\zeta^{\mathrm{left}}(L_{33})=Fa_{3}$, $[L_{33},L_{33}]=Fa_{1}\oplus Fa_{3}$, $\zeta^{\mathrm{right}}(L_{33})=\zeta(L_{33})=\langle0\rangle$. Moreover, if $\lambda a_{1}+\mu a_{2}+\nu a_{3}\not\in A$, then $\lambda^{2}+\lambda\nu-\lambda\mu\gamma+\mu^{2}\gamma+\mu\nu(2-\gamma)\neq0$, $L_{33}$ is non-nilpotent.
\item[\upshape(xi)] $\mathrm{\mathbf{Lei}}_{34}(3,F)=L_{34}$ is a sum of ideal $A=Fa_{1}\oplus Fa_{3}$ and subalgebra $B=Fa_{2}\oplus Fa_{3}$. Moreover, $A,B$ are non-nilpotent cyclic Leibniz algebras of dimension $2$, $[A,B]=[B,A]=Fa_{1}\oplus Fa_{3}$, $L_{34}=Fa_{1}\oplus Fa_{2}\oplus Fa_{3}$ where $[a_{1},a_{1}]=[a_{1},a_{3}]=a_{3}$, $[a_{1},a_{2}]=-a_{1}+\beta a_{3}\ (\beta=(\alpha+\gamma)(1+\alpha)^{-1})$, $[a_{2},a_{1}]=a_{1}+\alpha a_{3}\ (\alpha\neq0,\alpha\neq-1,\alpha\neq-2)$, $[a_{2},a_{2}]=\gamma a_{3}\ (\gamma\neq0)$, $[a_{2},a_{3}]=(2+\alpha)a_{3}$, $[a_{3},a_{1}]=[a_{3},a_{2}]=[a_{3},a_{3}]=0$, $\mathrm{\mathbf{Leib}}(L_{34})=\zeta^{\mathrm{left}}(L_{34})=Fa_{3}$, $[L_{34},L_{34}]=Fa_{1}\oplus Fa_{3}$, $\zeta^{\mathrm{right}}(L_{34})=\zeta(L_{34})=\langle0\rangle$. Moreover, if $\lambda a_{1}+\mu a_{2}+\nu a_{3}\not\in A$, then $\lambda^{2}+\lambda\mu\beta+\lambda\nu+\lambda\mu\alpha+\mu^{2}\gamma+\mu\nu(2+\alpha)\neq0$, $L_{34}$ is non-nilpotent.
\item[\upshape(xii)] $\mathrm{\mathbf{Lei}}_{35}(3,F)=L_{35}$ is a direct sum of ideal $B=Fa_{2}$ and a cyclic non-nilpotent subalgebra $A=Fa_{1}\oplus Fa_{3}$ of dimension $2$, $[A,B]=[B,A]=Fa_{2}$, so that $L_{35}=Fa_{1}\oplus Fa_{2}\oplus Fa_{3}$ where $[a_{1},a_{1}]=\gamma a_{3}\ (\gamma\neq0)$, $[a_{1},a_{2}]=a_{2}$, $[a_{1},a_{3}]=a_{3}$, $[a_{2},a_{1}]=-a_{2}$, $[a_{2},a_{2}]=[a_{2},a_{3}]=[a_{3},a_{1}]=[a_{3},a_{2}]=[a_{3},a_{3}]=0$, $\mathrm{\mathbf{Leib}}(L_{35})=\zeta^{\mathrm{left}}(L_{35})=Fa_{3}$, $[L_{35},L_{35}]=Fa_{2}\oplus Fa_{3}$, $\zeta^{\mathrm{right}}(L_{35})=\zeta(L_{35})=\langle0\rangle$, $L_{35}$ is non-nilpotent.
\item[\upshape(xiii)] $\mathrm{\mathbf{Lei}}_{36}(3,F)=L_{36}$ is a sum of abelian ideal $B=Fa_{2}\oplus Fa_{3}$ and a cyclic non-nilpotent subalgebra $A=Fa_{1}\oplus Fa_{3}$ of dimension $2$, $[A,B]=[B,A]=Fa_{2}\oplus Fa_{3}$, so that $L_{36}=Fa_{1}\oplus Fa_{2}\oplus Fa_{3}$ where $[a_{1},a_{1}]=\gamma a_{3}\ (\gamma\neq0)$, $[a_{1},a_{2}]=a_{2}$, $[a_{1},a_{3}]=a_{3}$, $[a_{2},a_{1}]=-a_{2}+\beta a_{3}\ (\beta\neq0)$, $[a_{2},a_{2}]=[a_{2},a_{3}]=[a_{3},a_{1}]=[a_{3},a_{2}]=[a_{3},a_{3}]=0$, $\mathrm{\mathbf{Leib}}(L_{36})=\zeta^{\mathrm{left}}(L_{36})=Fa_{3}$, $[L_{36},L_{36}]=Fa_{2}\oplus Fa_{3}$, $\zeta^{\mathrm{right}}(L_{36})=\zeta(L_{36})=\langle0\rangle$, $L_{36}$ is non-nilpotent.
\end{enumerate}
\end{thm}
\pf As in the previous theorem, we can see that $\mathrm{\mathbf{dim}}_{F}(\zeta(L))\leqslant1$. Suppose that $\mathrm{\mathbf{dim}}_{F}(\zeta(L))=1$. Since $L$ is not a Lie algebra, there is an element $a_{1}$ such that $[a_{1},a_{1}]=a_{3}\neq0$. Then, $\mathrm{\mathbf{Leib}}(L)=Fa_{3}$. Let $A=\langle a_{1}\rangle=Fa_{1}\oplus Fa_{3}$. An equality $\mathrm{\mathbf{Leib}}(L)\cap\zeta(L)=\langle0\rangle$ implies that $A\cap\zeta(L)=\langle0\rangle$, so that $L=A\oplus\zeta(L)$ and $A$ is an ideal of $L$. But, in this case, the factor-algebra $L/\mathrm{\mathbf{Leib}}(L)$ is abelian, and we obtain a contradiction. This contradiction shows that $\zeta(L)=\langle0\rangle$.

By what is noted above, $L/\mathrm{\mathbf{Leib}}(L)$ has an ideal $C/\mathrm{\mathbf{Leib}}(L)$ of dimension 1 (i.e., $C=Fc\oplus\mathrm{\mathbf{Leib}}(L)$ for some element $c$). If $[c,c]\neq0$ without loss of generality, we can put $c=a_{1}$. Then, subalgebra $\langle a_{1}\rangle=Fa_{1}\oplus Fa_{3}=A$ is an ideal which has codimension 1. Then, for every element $b$ such that $b\not\in A$, we have $L=A\oplus Fb$. By what is noted above, in this case, $[b,a_{1}]\in a_{1}+Fa_{3}$, so that $[b,a_{1}]=a_{1}+\alpha a_{3}$ for some element $\alpha\in F$. We have also $[a_{1},b]=-a_{1}+\beta a_{3}$ for some element $\beta\in F$. Since $[b,b]\in\mathrm{\mathbf{Leib}}(L)$, $[b,b]=\gamma a_{3}$ for some element $\gamma\in F$.

Suppose first that $\gamma=0$. In other words, we suppose that there exists an element $b$ such that $b\not\in A$ and $[b,b]=0$. Consider first the case when $A$ is nilpotent. We have:
\begin{gather*}
[b,a_{3}]=[b,[a_{1},a_{1}]]=[[b,a_{1}],a_{1}]+[a_{1},[b,a_{1}]]=\\
[a_{1}+\alpha a_{3},a_{1}]+[a_{1},a_{1}+\alpha a_{3}]=[a_{1},a_{1}]+[a_{1},a_{1}]=2a_{3}.
\end{gather*}
In particular, if we suppose that $\mathrm{\mathbf{char}}(F)=2$, then $[b,a_{3}]=0$. Since $[a_{1},a_{3}]=[a_{3},a_{1}]=[a_{3},b]=0$, $\mathrm{\mathbf{Leib}}(L)=Fa_{3}\leqslant\zeta(L)$, and we obtain a contradiction. This contradiction shows that $\mathrm{\mathbf{char}}(F)\neq2$. Further,
$$[b,[a_{1},b]]=[[b,a_{1}],b]+[a_{1},[b,b]]=[a_{1}+\alpha a_{3},b]=[a_{1},b]=-a_{1}+\beta a_{3}.$$
On the other hand,
\begin{gather*}
[b,[a_{1},b]]=[b,-a_{1}+\beta a_{3}]=-[b,a_{1}]+\beta[b,a_{3}]=\\
-(a_{1}+\alpha a_{3})+2\beta a_{3}=-a_{1}+(2\beta-\alpha)a_{3}.
\end{gather*}
It follows that $2\beta-\alpha=\beta$ or $\beta=\alpha$ and $[a_{1},b]=-a_{1}+\alpha a_{3}$.

Suppose that $\alpha=0$. Put $a_{2}=b$. Then $[a_{1},a_{2}]=-a_{1}$, $[a_{2},a_{1}]=a_{1}$, $[a_{2},a_{2}]=0$, $[a_{2},a_{3}]=2a_{3}$. Thus, we come to the following type of Leibniz algebras:
\begin{gather*}
L_{24}=Fa_{1}\oplus Fa_{2}\oplus Fa_{3}\ \mbox{where }[a_{1},a_{1}]=a_{3},\\
[a_{1},a_{2}]=-a_{1},[a_{2},a_{1}]=a_{1},[a_{2},a_{3}]=2a_{3},\\
[a_{1},a_{3}]=[a_{2},a_{2}]=[a_{3},a_{1}]=[a_{3},a_{2}]=[a_{3},a_{3}]=0.
\end{gather*}
Note also that $\mathrm{\mathbf{Leib}}(L_{24})=\zeta^{\mathrm{left}}(L_{24})=Fa_{3}$, $\zeta^{\mathrm{right}}(L_{24})=\zeta(L_{24})=\langle0\rangle$, $[L_{24},L_{24}]=Fa_{1}\oplus Fa_{3}$, $\mathrm{\mathbf{char}}(F)\neq2$, $L_{24}$ is non-nilpotent.

Suppose now that $\alpha\neq0$. Put again $a_{2}=b$. Thus, we come to the following type of Leibniz algebras:
\begin{gather*}
L_{25}=Fa_{1}\oplus Fa_{2}\oplus Fa_{3}\ \mbox{where }[a_{1},a_{1}]=a_{3},\\
[a_{1},a_{2}]=-a_{1}+\alpha a_{3},[a_{2},a_{1}]=a_{1}+\alpha a_{3}\ (\alpha\neq0),\\
[a_{2},a_{3}]=2a_{3},[a_{1},a_{3}]=[a_{2},a_{2}]=[a_{3},a_{1}]=[a_{3},a_{2}]=[a_{3},a_{3}]=0.
\end{gather*}
Note also that $\mathrm{\mathbf{Leib}}(L_{25})=\zeta^{\mathrm{left}}(L_{25})=Fa_{3}$, $\zeta^{\mathrm{right}}(L_{25})=\zeta(L_{25})=\langle0\rangle$, $[L_{25},L_{25}]=Fa_{1}\oplus Fa_{3}$, $\mathrm{\mathbf{char}}(F)\neq2$, $L_{25}$ is non-nilpotent.

Consider now the case when $A$ is not nilpotent. As we have seen above, we can choose an element $a_{1}$ such that $[a_{1},a_{3}]=a_{3}$. Then,
\begin{gather*}
[b,a_{3}]=[b,[a_{1},a_{1}]]=[[b,a_{1}],a_{1}]+[a_{1},[b,a_{1}]]=\\
[a_{1}+\alpha a_{3},a_{1}]+[a_{1},a_{1}+\alpha a_{3}]=\\
[a_{1},a_{1}]+[a_{1},a_{1}]+\alpha[a_{1},a_{3}]=(2+\alpha)a_{3}.
\end{gather*}
Further,
\begin{gather*}
[b,[a_{1},b]]=[[b,a_{1}],b]+[a_{1},[b,b]]=\\
[a_{1}+\alpha a_{3},b]=[a_{1},b]=-a_{1}+\beta a_{3}.
\end{gather*}
On the other hand,
\begin{gather*}
[b,[a_{1},b]]=[b,-a_{1}+\beta a_{3}]=-[b,a_{1}]+\beta[b,a_{3}]=\\
-(a_{1}+\alpha a_{3})+\beta(2+\alpha)a_{3}=-a_{1}+(2\beta+\beta\alpha-\alpha)a_{3}.
\end{gather*}
It follows that $2\beta+\beta\alpha-\alpha=\beta$ or $\beta+\beta\alpha-\alpha=0$. It follows that $\beta(1+\alpha)=\alpha$. We consider separately the case when $\alpha=0$. Then, $\beta=0$. If we suppose that $\mathrm{\mathbf{char}}(F)=2$, then $[b,a_{3}]=0$. Since $[a_{3},b]=0$, $b\in\mathrm{\mathbf{Ann}}_{L}(Fa_{3})=\mathrm{\mathbf{Ann}}_{L}(\mathrm{\mathbf{Leib}}(L))$. But, $\mathrm{\mathbf{Ann}}_{L}(\mathrm{\mathbf{Leib}}(L))$ is an ideal of $L$. Then, $[b,a_{1}]=a_{1}\in Fb\oplus Fa_{3}$, and we obtain a contradiction. This contradiction shows that $\mathrm{\mathbf{char}}(F)\neq2$. Put $a_{2}=b$. We come to the following type of Leibniz algebras:
\begin{gather*}
L_{26}=Fa_{1}\oplus Fa_{2}\oplus Fa_{3}\ \mbox{where }[a_{1},a_{1}]=[a_{1},a_{3}]=a_{3},[a_{1},a_{2}]=-a_{1},\\
[a_{2},a_{1}]=a_{1},[a_{2},a_{3}]=2a_{3},[a_{2},a_{2}]=[a_{3},a_{1}]=[a_{3},a_{2}]=[a_{3},a_{3}]=0.
\end{gather*}
Note also that $\mathrm{\mathbf{Leib}}(L_{26})=\zeta^{\mathrm{left}}(L_{26})=Fa_{3}$, $[L_{26},L_{26}]=Fa_{1}\oplus Fa_{3}$, $\zeta^{\mathrm{right}}(L_{26})=\zeta(L_{26})=\langle0\rangle$, $\mathrm{\mathbf{char}}(F)\neq2$, $L_{26}$ is non-nilpotent.

Suppose that $\alpha\neq0$. The equality $\beta(1+\alpha)=\alpha$ shows that $\alpha\neq-1$. In this case, $\beta=\alpha(1+\alpha)^{-1}$. If we suppose that $\alpha=-2$, then $[b,a_{3}]=0$. Using the above arguments, we obtain a contradiction. This contradiction shows that $\alpha\neq-2$. Put $a_{2}=b$. We come to the following type of Leibniz algebras:
\begin{gather*}
L_{27}=Fa_{1}\oplus Fa_{2}\oplus Fa_{3}\ \mbox{where }[a_{1},a_{1}]=[a_{1},a_{3}]=a_{3},\\
[a_{1},a_{2}]=-a_{1}+\beta a_{3}\ (\beta=\alpha(1+\alpha)^{-1}),[a_{2},a_{1}]=a_{1}+\alpha a_{3},\\
[a_{2},a_{3}]=(2+\alpha)a_{3}\ (\alpha\neq0,\alpha\neq-1,\alpha\neq-2),\\
[a_{2},a_{2}]=[a_{3},a_{1}]=[a_{3},a_{2}]=[a_{3},a_{3}]=0.
\end{gather*}
Note also that $\mathrm{\mathbf{Leib}}(L_{27})=\zeta^{\mathrm{left}}(L_{27})=Fa_{3}$, $[L_{27},L_{27}]=Fa_{1}\oplus Fa_{3}$, $\zeta^{\mathrm{right}}(L_{27})=\zeta(L_{27})=\langle0\rangle$, $L_{27}$ is non-nilpotent.

Suppose now that $[b,b]\neq0$ for every element $b$ such that $b\not\in A$, so that $[b,b]=\gamma a_{3}$ and $\gamma\neq0$. Consider first the case when $A$ is nilpotent. We have:
\begin{gather*}
[b,a_{3}]=[b,[a_{1},a_{1}]]=[[b,a_{1}],a_{1}]+[a_{1},[b,a_{1}]]=\\
[a_{1}+\alpha a_{3},a_{1}]+[a_{1},a_{1}+\alpha a_{3}]=[a_{1},a_{1}]+[a_{1},a_{1}]=2a_{3}.
\end{gather*}
Again, we obtain that $\mathrm{\mathbf{char}}(F)\neq2$. Using the arguments above, we can get that $[a_{1},b]=-a_{1}+\alpha a_{3}$ again. We consider separately the case when $\alpha=0$ (that is, $[a_{1},b]=-a_{1}$, $[b,a_{1}]=a_{1}$). Let $\lambda a_{1}+\mu b+\nu a_{3}$ be the arbitrary element of $L$. We have:
\begin{gather*}
[\lambda a_{1}+\mu b+\nu a_{3},\lambda a_{1}+\mu b+\nu a_{3}]=\\
\lambda^{2}[a_{1},a_{1}]+\lambda\mu[a_{1},b]+\lambda\mu[b,a_{1}]+\mu^{2}[b,b]+\mu\nu[b,a_{3}]=\\
\lambda^{2}a_{3}-\lambda\mu a_{1}+\lambda\mu a_{1}+\mu^{2}\gamma a_{3}+2\mu\nu a_{3}=\\
(\lambda^{2}+\mu^{2}\gamma+2\mu\nu)a_{3}.
\end{gather*}
Then, a condition $[b,b]\neq0$ for every element $b\not\in A$ yields that $\lambda^{2}+\mu^{2}\gamma+2\mu\nu\neq0$. Put $a_{2}=b$. Thus, we come to the following type of Leibniz algebras:
\begin{gather*}
L_{28}=Fa_{1}\oplus Fa_{2}\oplus Fa_{3}\ \mbox{where }[a_{1},a_{1}]=a_{3},[a_{1},a_{2}]=-a_{1},[a_{2},a_{1}]=a_{1},\\
[a_{2},a_{2}]=\gamma a_{3}\ (\gamma\neq0),[a_{2},a_{3}]=2a_{3},\\
[a_{1},a_{3}]=[a_{3},a_{1}]=[a_{3},a_{2}]=[a_{3},a_{3}]=0.
\end{gather*}
Note also that $\mathrm{\mathbf{Leib}}(L_{28})=\zeta^{\mathrm{left}}(L_{28})=Fa_{3}$, $[L_{28},L_{28}]=Fa_{1}\oplus Fa_{3}$, $\zeta^{\mathrm{right}}(L_{28})=\zeta(L_{28})=\langle0\rangle$, $\mathrm{\mathbf{char}}(F)\neq2$. Moreover, if $\lambda a_{1}+\mu a_{2}+\nu a_{3}\not\in A$, then $\lambda^{2}+\mu^{2}\gamma+2\mu\nu\neq0$, $L_{28}$ is non-nilpotent.

Suppose that $\alpha\neq0$. Let $\lambda a_{1}+\mu b+\nu a_{3}$ be the arbitrary element of $L$. We have:
\begin{gather*}
[\lambda a_{1}+\mu b+\nu a_{3},\lambda a_{1}+\mu b+\nu a_{3}]=\\
\lambda^{2}[a_{1},a_{1}]+\lambda\mu[a_{1},b]+\lambda\mu[b,a_{1}]+\mu^{2}[b,b]+\mu\nu[b,a_{3}]=\\
\lambda^{2}a_{3}+\lambda\mu(-a_{1}+\alpha a_{3})+\lambda\mu(a_{1}+\alpha a_{3})+\mu^{2}\gamma a_{3}+2\mu\nu a_{3}=\\
\lambda^{2}a_{3}-\lambda\mu a_{1}+\lambda\mu\alpha a_{3}+\lambda\mu a_{1}+\lambda\mu\alpha a_{3}+\mu^{2}\gamma a_{3}+2\mu\nu a_{3}=\\
(\lambda^{2}+2\lambda\mu\alpha+\mu^{2}\gamma+2\mu\nu)a_{3}.
\end{gather*}
As above, $\lambda^{2}+2\lambda\mu\alpha+\mu^{2}\gamma+2\mu\nu\neq0$. Put again $a_{2}=b$. Then, we come to the following type of Leibniz algebras:
\begin{gather*}
L_{29}=Fa_{1}\oplus Fa_{2}\oplus Fa_{3}\ \mbox{where }[a_{1},a_{1}]=a_{3},\\
[a_{1},a_{2}]=-a_{1}+\alpha a_{3},[a_{2},a_{1}]=a_{1}+\alpha a_{3}\ (\alpha\neq0),\\
[a_{2},a_{2}]=\gamma a_{3}\ (\gamma\neq0),[a_{2},a_{3}]=2a_{3},\\
[a_{1},a_{3}]=[a_{3},a_{1}]=[a_{3},a_{2}]=[a_{3},a_{3}]=0.
\end{gather*}
Note also that $\mathrm{\mathbf{Leib}}(L_{29})=\zeta^{\mathrm{left}}(L_{29})=Fa_{3}$, $[L_{29},L_{29}]=Fa_{1}\oplus Fa_{3}$, $\zeta^{\mathrm{right}}(L_{29})=\zeta(L_{29})=\langle0\rangle$, $\mathrm{\mathbf{char}}(F)\neq2$. Moreover, if $\lambda a_{1}+\mu a_{2}+\nu a_{3}\not\in A$, then $\lambda^{2}+2\lambda\mu\alpha+\mu^{2}\gamma+2\mu\nu\neq0$, $L_{29}$ is non-nilpotent.

Suppose, now, that $[b,b]\neq0$ for every element $b$ such that $b\not\in A$ and a subalgebra $A$ is not nilpotent. As we have seen above, we can choose an element $a_{1}$ such that $[a_{1},a_{3}]=a_{3}$. Then,
\begin{gather*}
[b,a_{3}]=[b,[a_{1},a_{1}]]=[[b,a_{1}],a_{1}]+[a_{1},[b,a_{1}]]=\\
[a_{1}+\alpha a_{3},a_{1}]+[a_{1},a_{1}+\alpha a_{3}]=\\
[a_{1},a_{1}]+[a_{1},a_{1}]+\alpha[a_{1},a_{3}]=(2+\alpha)a_{3}.
\end{gather*}
If we suppose that $\alpha=-2$, then $[b,a_{3}]=0$. Since $[a_{3},b]=0$, $b\in\mathrm{\mathbf{Ann}}_{L}(Fa_{3})=\mathrm{\mathbf{Ann}}_{L}(\mathrm{\mathbf{Leib}}(L))$. But, $\mathrm{\mathbf{Ann}}_{L}(\mathrm{\mathbf{Leib}}(L))$ is an ideal of $L$. Then, $[b,a_{1}]\in Fb\oplus Fa_{3}$, and we obtain a contradiction. This contradiction shows that $\alpha\neq-2$. Further,
\begin{gather*}
[b,[a_{1},b]]=[[b,a_{1}],b]+[a_{1},[b,b]]=[a_{1}+\alpha a_{3},b]+[a_{1},\gamma a_{3}]=\\
[a_{1},b]+\gamma[a_{1},a_{3}]=-a_{1}+\beta a_{3}+\gamma a_{3}=-a_{1}+(\beta+\gamma)a_{3}.
\end{gather*}
On the other hand,
\begin{gather*}
[b,[a_{1},b]]=[b,-a_{1}+\beta a_{3}]=-[b,a_{1}]+\beta[b,a_{3}]=\\
-(a_{1}+\alpha a_{3})+\beta(2+\alpha)a_{3}=-a_{1}+(2\beta+\beta\alpha-\alpha)a_{3}.
\end{gather*}
It follows that $2\beta+\beta\alpha-\alpha=\beta+\gamma$ or $\beta(1+\alpha)=\alpha+\gamma$. We consider separately the case when $\alpha=0$ and $\alpha=-1$.

Let $\alpha=0$, then $\beta=\gamma$ and $[b,a_{3}]=2a_{3}$. Using the arguments given above, we can obtain that $\mathrm{\mathbf{char}}(F)\neq2$. In this case, $[b,a_{1}]=a_{1}$ and $[a_{1},b]=-a_{1}+\gamma a_{3}$. Let $\lambda a_{1}+\mu b+\nu a_{3}$ be the arbitrary element of $L$. We have:
\begin{gather*}
[\lambda a_{1}+\mu b+\nu a_{3},\lambda a_{1}+\mu b+\nu a_{3}]=\\
\lambda^{2}[a_{1},a_{1}]+\lambda\mu[a_{1},b]+\lambda\nu[a_{1},a_{3}]+\lambda\mu[b,a_{1}]+\mu^{2}[b,b]+\mu\nu[b,a_{3}]=\\
\lambda^{2}a_{3}+\lambda\mu(-a_{1}+\gamma a_{3})+\lambda\mu a_{3}+\lambda\mu a_{1}+\mu^{2}\gamma a_{3}+2\mu\nu a_{3}=\\
\lambda^{2}a_{3}-\lambda\mu a_{1}+\lambda\mu\gamma a_{3}+\lambda\nu a_{3}+\lambda\mu a_{1}+\mu^{2}\gamma a_{3}+2\mu\nu a_{3}=\\
(\lambda^{2}+\lambda\mu\gamma+\lambda\nu+\mu^{2}\gamma+2\mu\nu)a_{3}.
\end{gather*}
Then, the condition $[b,b]\neq0$ for every element $b\not\in A$ yields that $\lambda^{2}+\lambda\mu\gamma+\lambda\nu+\mu^{2}\gamma+2\mu\nu\neq0$. Put $a_{2}=b$. Thus, we come to the following type of Leibniz algebras:
\begin{gather*}
L_{30}=Fa_{1}\oplus Fa_{2}\oplus Fa_{3}\ \mbox{where }[a_{1},a_{1}]=[a_{1},a_{3}]=a_{3},\\
[a_{1},a_{2}]=-a_{1}+\gamma a_{3}\ (\gamma\neq0),[a_{2},a_{1}]=a_{1},[a_{2},a_{2}]=\gamma a_{3},[a_{2},a_{3}]=2a_{3},\\
[a_{3},a_{1}]=[a_{3},a_{2}]=[a_{3},a_{3}]=0.
\end{gather*}
Note also that $\mathrm{\mathbf{Leib}}(L_{30})=\zeta^{\mathrm{left}}(L_{30})=Fa_{3}$, $[L_{30},L_{30}]=Fa_{1}\oplus Fa_{3}$, $\zeta^{\mathrm{right}}(L_{30})=\zeta(L_{30})=\langle0\rangle$, $\mathrm{\mathbf{char}}(F)\neq2$. Moreover, if $\lambda a_{1}+\mu a_{2}+\nu a_{3}\not\in A$, then $\lambda^{2}+\lambda\mu\gamma+\lambda\nu+\mu^{2}\gamma+2\mu\nu\neq0$, $L_{30}$ is non-nilpotent.

Let now $\alpha=-1$. Then $\gamma=1$, $[b,b]=[b,a_{3}]=a_{3}$, $[b,a_{1}]=a_{1}-a_{3}$, $[a_{1},b]=-a_{1}+\beta a_{3}$. If $\beta=0$, then $[a_{1},b]=-a_{1}$. Let $\lambda a_{1}+\mu b+\nu a_{3}$ be the arbitrary element of $L$. We have:
\begin{gather*}
[\lambda a_{1}+\mu b+\nu a_{3},\lambda a_{1}+\mu b+\nu a_{3}]=\\
\lambda^{2}[a_{1},a_{1}]+\lambda\mu[a_{1},b]+\lambda\nu[a_{1},a_{3}]+\lambda\mu[b,a_{1}]+\mu^{2}[b,b]+\mu\nu[b,a_{3}]=\\
\lambda^{2}a_{3}-\lambda\mu a_{1}+\lambda\nu a_{3}+\lambda\mu(a_{1}-a_{3})+\mu^{2}a_{3}+\mu\nu a_{3}=\\
\lambda^{2}a_{3}-\lambda\mu a_{1}+\lambda\nu a_{3}+\lambda\mu a_{1}-\lambda\mu a_{3}+\mu^{2}a_{3}+\mu\nu a_{3}=\\
(\lambda^{2}+\lambda\nu-\lambda\mu+\mu^{2}+\mu\nu)a_{3}.
\end{gather*}
As above, $\lambda^{2}+\lambda\nu-\lambda\mu+\mu^{2}+\mu\nu\neq0$. Put $a_{2}=b$. Thus, we come to the following type of Leibniz algebras:
\begin{gather*}
L_{31}=Fa_{1}\oplus Fa_{2}\oplus Fa_{3}\ \mbox{where}\\
[a_{1},a_{1}]=[a_{1},a_{3}]=[a_{2},a_{2}]=[a_{2},a_{3}]=a_{3},\\
[a_{1},a_{2}]=-a_{1},[a_{2},a_{1}]=a_{1}-a_{3},\\
[a_{3},a_{1}]=[a_{3},a_{2}]=[a_{3},a_{3}]=0.
\end{gather*}
Note also that $\mathrm{\mathbf{Leib}}(L_{31})=\zeta^{\mathrm{left}}(L_{31})=Fa_{3}$, $[L_{31},L_{31}]=Fa_{1}\oplus Fa_{3}$, $\zeta^{\mathrm{right}}(L_{31})=\zeta(L_{31})=\langle0\rangle$. Moreover, if $\lambda a_{1}+\mu a_{2}+\nu a_{3}\not\in A$, then $\lambda^{2}+\lambda\nu-\lambda\mu+\mu^{2}+\mu\nu\neq0$, $L_{31}$ is non-nilpotent.

Suppose that $\beta\neq0$. Let $\lambda a_{1}+\mu b+\nu a_{3}$ be the arbitrary element of $L$. We have:
\begin{gather*}
[\lambda a_{1}+\mu b+\nu a_{3},\lambda a_{1}+\mu b+\nu a_{3}]=\\
\lambda^{2}[a_{1},a_{1}]+\lambda\mu[a_{1},b]+\lambda\nu[a_{1},a_{3}]+\lambda\mu[b,a_{1}]+\mu^{2}[b,b]+\mu\nu[b,a_{3}]=\\
\lambda^{2}a_{3}+\lambda\mu(-a_{1}+\beta a_{3})+\lambda\nu a_{3}+\lambda\mu(a_{1}-a_{3})+\mu^{2}a_{3}+\mu\nu a_{3}=\\
\lambda^{2}a_{3}-\lambda\mu a_{1}+\lambda\mu\beta a_{3}+\lambda\nu a_{3}+\lambda\mu a_{1}-\lambda\mu a_{3}+\mu^{2}a_{3}+\mu\nu a_{3}=\\
(\lambda^{2}+\lambda\mu\beta+\lambda\nu-\lambda\mu+\mu^{2}+\mu\nu)a_{3}.
\end{gather*}
As above, $\lambda^{2}+\lambda\mu\beta+\lambda\nu-\lambda\mu+\mu^{2}+\mu\nu\neq0$. Put $a_{2}=b$. Thus, we come to the following type of Leibniz algebras:
\begin{gather*}
L_{32}=Fa_{1}\oplus Fa_{2}\oplus Fa_{3}\ \mbox{where}\\
[a_{1},a_{1}]=[a_{1},a_{3}]=[a_{2},a_{2}]=[a_{2},a_{3}]=a_{3},\\
[a_{1},a_{2}]=-a_{1}+\beta a_{3}\ (\beta\neq0),[a_{2},a_{1}]=a_{1}-a_{3},\\
[a_{3},a_{1}]=[a_{3},a_{2}]=[a_{3},a_{3}]=0.
\end{gather*}
Note also that $\mathrm{\mathbf{Leib}}(L_{32})=\zeta^{\mathrm{left}}(L_{32})=Fa_{3}$, $[L_{32},L_{32}]=Fa_{1}\oplus Fa_{3}$, $\zeta^{\mathrm{right}}(L_{32})=\zeta(L_{32})=\langle0\rangle$. Moreover, if $\lambda a_{1}+\mu a_{2}+\nu a_{3}\not\in A$, then $\lambda^{2}+\lambda\mu\beta+\lambda\nu-\lambda\mu+\mu^{2}+\mu\nu\neq0$, $L_{32}$ is non-nilpotent.

Suppose now that $\alpha\neq0$ and $\alpha\neq-1$. As we have noted above, $\alpha\neq-2$. We obtain $\beta=(\alpha+\gamma)(1+\alpha)^{-1}$. If $\beta=0$, then $\alpha=-\gamma$, $[a_{1},b]=-a_{1}$, $[b,a_{1}]=a_{1}-\gamma a_{3}$, $[b,a_{3}]=(2-\gamma)a_{3}$. Let $\lambda a_{1}+\mu b+\nu a_{3}$ be the arbitrary element of $L$. We have:
\begin{gather*}
[\lambda a_{1}+\mu b+\nu a_{3},\lambda a_{1}+\mu b+\nu a_{3}]=\\
\lambda^{2}[a_{1},a_{1}]+\lambda\mu[a_{1},b]+\lambda\nu[a_{1},a_{3}]+\lambda\mu[b,a_{1}]+\mu^{2}[b,b]+\mu\nu[b,a_{3}]=\\
\lambda^{2}a_{3}-\lambda\mu a_{1}+\lambda\nu a_{3}+\lambda\mu(a_{1}-\gamma a_{3})+\mu^{2}\gamma a_{3}+\mu\nu(2-\gamma)a_{3}=\\
\lambda^{2}a_{3}-\lambda\mu a_{1}+\lambda\nu a_{3}+\lambda\mu a_{1}-\lambda\mu\gamma a_{3}+\mu^{2}\gamma a_{3}+\mu\nu(2-\gamma)a_{3}=\\ (\lambda^{2}+\lambda\nu-\lambda\mu\gamma+\mu^{2}\gamma+\mu\nu(2-\gamma))a_{3}.
\end{gather*}
As above, $\lambda^{2}+\lambda\nu-\lambda\mu\gamma+\mu^{2}\gamma+\mu\nu(2-\gamma)\neq0$. Put $a_{2}=b$. Thus, we come to the following type of Leibniz algebras:
\begin{gather*}
L_{33}=Fa_{1}\oplus Fa_{2}\oplus Fa_{3}\ \mbox{where }[a_{1},a_{1}]=[a_{1},a_{3}]=a_{3},\\
[a_{1},a_{2}]=-a_{1},[a_{2},a_{1}]=a_{1}-\gamma a_{3}\ (\gamma\neq0,\gamma\neq1,\gamma\neq2),\\
[a_{2},a_{2}]=\gamma a_{3},[a_{2},a_{3}]=(2-\gamma)a_{3},[a_{3},a_{1}]=[a_{3},a_{2}]=[a_{3},a_{3}]=0.
\end{gather*}
Note also that $\mathrm{\mathbf{Leib}}(L_{33})=\zeta^{\mathrm{left}}(L_{33})=Fa_{3}$, $[L_{33},L_{33}]=Fa_{1}\oplus Fa_{3}$, $\zeta^{\mathrm{right}}(L_{33})=\zeta(L_{33})=\langle0\rangle$. Moreover, if $\lambda a_{1}+\mu a_{2}+\nu a_{3}\not\in A$, then $\lambda^{2}+\lambda\nu-\lambda\mu\gamma+\mu^{2}\gamma+\mu\nu(2-\gamma)\neq0$, $L_{33}$ is non-nilpotent.

Suppose that $\beta\neq0$. Then $[a_{1},b]=-a_{1}+\beta a_{3}\ (\beta=(\alpha+\gamma)(1+\alpha)^{-1})$, $[b,a_{1}]=a_{1}+\alpha a_{3}$, $[b,a_{3}]=(2+\alpha)a_{3}$. Let $\lambda a_{1}+\mu b+\nu a_{3}$ be the arbitrary element of $L$. We have:
\begin{gather*}
[\lambda a_{1}+\mu b+\nu a_{3},\lambda a_{1}+\mu b+\nu a_{3}]=\\
\lambda^{2}[a_{1},a_{1}]+\lambda\mu[a_{1},b]+\lambda\nu[a_{1},a_{3}]+\lambda\mu[b,a_{1}]+\mu^{2}[b,b]+\mu\nu[b,a_{3}]=\\
\lambda^{2}a_{3}+\lambda\mu(-a_{1}+\beta a_{3})+\lambda\nu a_{3}+\lambda\mu(a_{1}+\alpha a_{3})+\mu^{2}\gamma a_{3}+\mu\nu(2+\alpha)a_{3}=\\
\lambda^{2}a_{3}-\lambda\mu a_{1}+\lambda\mu\beta a_{3}+\lambda\nu a_{3}+\lambda\mu a_{1}+\lambda\mu\alpha a_{3}+\mu^{2}\gamma a_{3}+\mu\nu(2+\alpha)a_{3}=\\
(\lambda^{2}+\lambda\mu\beta+\lambda\nu+\lambda\mu\alpha+\mu^{2}\gamma+\mu\nu(2+\alpha))a_{3}.
\end{gather*}
As above, $\lambda^{2}+\lambda\mu\beta+\lambda\nu+\lambda\mu\alpha+\mu^{2}\gamma+\mu\nu(2+\alpha)\neq0$. Put $a_{2}=b$. Thus, we come to the following type of Leibniz algebras:
\begin{gather*}
L_{34}=Fa_{1}\oplus Fa_{2}\oplus Fa_{3}\ \mbox{where }[a_{1},a_{1}]=[a_{1},a_{3}]=a_{3},\\
[a_{1},a_{2}]=-a_{1}+\beta a_{3}\ (\beta=(\alpha+\gamma)(1+\alpha)^{-1}),\\
[a_{2},a_{1}]=a_{1}+\alpha a_{3}\ (\alpha\neq0,\alpha\neq-1,\alpha\neq-2),\\
[a_{2},a_{2}]=\gamma a_{3}\ (\gamma\neq0),[a_{2},a_{3}]=(2+\alpha)a_{3},\\
[a_{3},a_{1}]=[a_{3},a_{2}]=[a_{3},a_{3}]=0.
\end{gather*}
Note also that $\mathrm{\mathbf{Leib}}(L_{34})=\zeta^{\mathrm{left}}(L_{34})=Fa_{3}$, $[L_{34},L_{34}]=Fa_{1}\oplus Fa_{3}$, $\zeta^{\mathrm{right}}(L_{34})=\zeta(L_{34})=\langle0\rangle$. Moreover, if $\lambda a_{1}+\mu a_{2}+\nu a_{3}\not\in A$, then $\lambda^{2}+\lambda\mu\beta+\lambda\nu+\lambda\mu\alpha+\mu^{2}\gamma+\mu\nu(2+\alpha)\neq0$, $L_{34}$ is non-nilpotent.

Suppose now that $[c,c]=0$. Put again $\mathrm{\mathbf{Leib}}(L)=Fa_{3}$. Let $b$ be an element such that $b\not\in C$.

Suppose that a subalgebra $C=Fc\oplus\mathrm{\mathbf{Leib}}(L)$ is not abelian. As we have seen above, we can choose an element $a_{3}$ such that $[c,a_{3}]=a_{3}$. Since $\mathrm{\mathbf{Leib}}(L)=Fa_{3}$ is an ideal, $[b,a_{3}]=\eta a_{3}$ for some element $\eta\in F$. Using the above arguments without loss of generality we can assume that $[b,c]\in c+Fa_{3}$, $[c,b]\in-c+Fa_{3}$, so that $[b,c]=c+\alpha a_{3}$, $[c,b]=-c+\beta a_{3}$ for some elements $\alpha,\beta\in F$. Then
\begin{gather*}
a_{3}=[c,a_{3}]=[c+\alpha a_{3},a_{3}]=[[b,c],a_{3}]=[b,[c,a_{3}]]-[c,[b,a_{3}]]=\\
[b,a_{3}]-[c,\eta a_{3}]=\eta a_{3}-\eta[c,a_{3}]=\eta a_{3}-\eta a_{3}=0.
\end{gather*}
This contradiction shows that a subalgebra $C$ is abelian.

Note that $\eta\neq0$. In fact, otherwise $a_{3}\in\zeta(L)$, and we obtain a contradiction. This contradiction shows that $[b,a_{3}]\neq0$ for every element $b\not\in C$. As we have seen above, we can choose an element $b$ such that $[b,a_{3}]=a_{3}$.

Suppose first that there exists an element $b\not\in C$ such that $[b,b]\neq0$. It follows that $[b,b]=\gamma a_{3}$ where $\gamma$ is a non-zero element of $F$. Then, subalgebra $B=Fb\oplus Fa_{3}$ is non-abelian. We have now
\begin{gather*}
[c,a_{3}]=[c,\gamma^{-1}[b,b]]=\gamma^{-1}[c,[b,b]]=\gamma^{-1}([[c,b],b]+[b,[c,b]])=\\
\gamma^{-1}([-c+\beta a_{3},b]+[b,-c+\beta a_{3}])=\gamma^{-1}(-[c,b]-[b,c]+\beta[b,a_{3}])=\\
\gamma^{-1}(c-\beta a_{3}-c-\alpha a_{3}+\beta a_{3})=-\alpha\gamma^{-1}a_{3}.
\end{gather*}
On the other hand, we proved above that $[c,a_{3}]=0$. Since $\gamma\neq0$, it follows that $\alpha=0$. Hence, $[b,c]=c$.

If $\beta=0$, then put $a_{1}=b$, $a_{2}=c$. A subalgebra $Fc$ is an ideal of $L$ and we come to the following type of Leibniz algebras:
\begin{gather*}
L_{35}=Fa_{1}\oplus Fa_{2}\oplus Fa_{3}\ \mbox{where }[a_{1},a_{1}]=\gamma a_{3}\ (\gamma\neq0),\\
[a_{1},a_{2}]=a_{2},[a_{1},a_{3}]=a_{3},[a_{2},a_{1}]=-a_{2},\\
[a_{2},a_{2}]=[a_{2},a_{3}]=[a_{3},a_{1}]=[a_{3},a_{2}]=[a_{3},a_{3}]=0.
\end{gather*}
Note also that $\mathrm{\mathbf{Leib}}(L_{35})=\zeta^{\mathrm{left}}(L_{35})=Fa_{3}$, $[L_{35},L_{35}]=Fa_{2}\oplus Fa_{3}$, $\zeta^{\mathrm{right}}(L_{35})=\zeta(L_{35})=\langle0\rangle$, $L_{35}$ is non-nilpotent.

If $\beta\neq0$, then put $a_{1}=b$, $a_{2}=c$. A subalgebra $Fc$ is not an ideal of $L$ and we come to the following type of Leibniz algebras:
\begin{gather*}
L_{36}=Fa_{1}\oplus Fa_{2}\oplus Fa_{3}\ \mbox{where }[a_{1},a_{1}]=\gamma a_{3}\ (\gamma\neq0),\\
[a_{1},a_{2}]=a_{2},[a_{1},a_{3}]=a_{3},[a_{2},a_{1}]=-a_{2}+\beta a_{3}\ (\beta\neq0),\\
[a_{2},a_{2}]=[a_{2},a_{3}]=[a_{3},a_{1}]=[a_{3},a_{2}]=[a_{3},a_{3}]=0.
\end{gather*}
Note also that $\mathrm{\mathbf{Leib}}(L_{36})=\zeta^{\mathrm{left}}(L_{36})=Fa_{3}$, $[L_{36},L_{36}]=Fa_{2}\oplus Fa_{3}$, $\zeta^{\mathrm{right}}(L_{36})=\zeta(L_{36})=\langle0\rangle$, $L_{36}$ is non-nilpotent.

Suppose now that $[b,b]=0$ for each element $b\not\in C$. Let $u=\lambda c+\mu b+\nu a_{3}$ be an arbitrary element of $L$. Then,
\begin{gather*}
[\lambda c+\mu b+\nu a_{3},\lambda c+\mu b+\nu a_{3}]=\lambda\mu[c,b]+\lambda\mu[b,c]+\mu\nu[b,a_{3}]=\\
\lambda\mu(-c+\beta a_{3})+\lambda\mu(c+\alpha a_{3})+\mu\nu a_{3}=\\
(\lambda\mu\beta+\lambda\mu\alpha+\mu\nu)a_{3}.
\end{gather*}
If $u\not\in C$, then $\mu\neq0$. If $\lambda=0$, $\mu=\nu=1$, then $[u,u]=a_{3}\neq0$, and we obtain a contradiction. This contradiction shows that this situation is not possible.
\qed

The following natural situation appears when $\mathrm{\mathbf{dim}}_{F}(\mathrm{\mathbf{Leib}}(L))=2$. Immediately, we obtain the following two subcases:

(IIA) the intersection $\zeta(L)\cap\mathrm{\mathbf{Leib}}(L)$ is not trivial;

(IIB) $\zeta(L)\cap\mathrm{\mathbf{Leib}}(L)=\langle0\rangle$.

\begin{thm}\label{T5}
Let $L$ be a nilpotent Leibniz algebra over a field $F$ having dimension $3$, which is not a Lie algebra. Suppose that $\mathrm{\mathbf{dim}}_{F}(\mathrm{\mathbf{Leib}}(L))=2$. Then, $L$ is an algebra of the following type:
\begin{enumerate}
\item[\upshape] $\mathrm{\mathbf{Lei}}_{37}(3,F)=L_{37}$ is a cyclic nilpotent Leibniz algebra, so that $L_{37}=Fa_{1}\oplus Fa_{2}\oplus Fa_{3}$ where $[a_{1},a_{1}]=a_{2}$, $[a_{1},a_{2}]=a_{3}$, $[a_{1},a_{3}]=[a_{2},a_{1}]=[a_{2},a_{2}]=[a_{2},a_{3}]=[a_{3},a_{1}]=[a_{3},a_{2}]=[a_{3},a_{3}]=0$, $\mathrm{\mathbf{Leib}}(L_{37})=[L_{37},L_{37}]=\zeta^{\mathrm{left}}(L_{37})=Fa_{2}\oplus Fa_{3}$, $\zeta^{\mathrm{right}}(L_{37})=\zeta(L_{37})=Fa_{3}$, $\mathrm{\mathbf{ncl}}(L_{37})=3$.
\end{enumerate}
\end{thm}
\pf Since $\mathrm{\mathbf{dim}}_{F}(L)=3$, $\mathrm{\mathbf{ncl}}(L)\leqslant3$. Suppose first that $\mathrm{\mathbf{ncl}}(L)=3$. Then, $L$ has an upper central series of a length 3:
$$\langle0\rangle=C_{0}\leqslant C_{1}\leqslant C_{2}\leqslant C_{3}=L.$$
Every factor of this series must be non-trivial. Therefore every factor of this series has dimension 1. Let $a_{1}$ be an element of $L$ such that $a_{1}\not\in C_{2}$. The fact that $L/C_{2}$ is abelian implies that $a_{2}=[a_{1},a_{1}]\in C_{2}$. Suppose that $a_{2}\in C_{1}$. Since $\mathrm{\mathbf{dim}}_{F}(C_{1})=1$, $C_{1}=Fa_{2}$. Choose an element $b\in C_{2}$ such that $b\not\in C_{1}$. Then, $[a_{1},b],[b,a_{1}],[b,b]\in C_{1}$, so that $[a_{1},b]=\alpha a_{2}$, $[b,a_{1}]=\beta a_{2}$, $[b,b]=\gamma a_{2}$ for some elements $\alpha,\beta,\gamma\in F$. It is not hard to see that the elements $a_{1},b,a_{2}$ generate $L$. Let $\lambda_{1}a_{1}+\lambda_{2}b+\lambda_{3}a_{2}$, $\mu_{1}a_{1}+\mu_{2}b+\mu_{3}a_{2}$ be two arbitrary elements of $L$. We have:
\begin{gather*}
[\lambda_{1}a_{1}+\lambda_{2}b+\lambda_{3}a_{2},\mu_{1}a_{1}+\mu_{2}b+\mu_{3}a_{2}]=\\
\lambda_{1}\mu_{1}[a_{1},a_{1}]+\lambda_{1}\mu_{2}[a_{1},b]+\lambda_{1}\mu_{3}[a_{1},a_{2}]+\\
\lambda_{2}\mu_{1}[b,a_{1}]+\lambda_{2}\mu_{2}[b,b]+\lambda_{2}\mu_{3}[b,a_{2}]=\\
\lambda_{1}\mu_{1}a_{2}+\lambda_{1}\mu_{2}\alpha a_{2}+\lambda_{2}\mu_{1}\beta a_{2}+\lambda_{2}\mu_{2}\gamma a_{2}\in C_{1}.
\end{gather*}
It follows that $[L,L]\leqslant C_{1}=\zeta(L)$ and hence, $\mathrm{\mathbf{ncl}}(L)=2$. This contradiction shows that $a_{2}\\not\in C_{1}$. Then, $a_{3}=[a_{1},a_{2}]\neq0$. It follows that $C_{2}=Fa_{2}\oplus Fa_{3}$. We have $[a_{2},a_{3}]=[a_{3},a_{2}]=0$, so that $C_{2}$ is an abelian subalgebra, $[a_{1},a_{3}]=[a_{3},a_{1}]=[a_{2},a_{1}]=0$. It follows that
\begin{gather*}
[a_{1}+a_{2},a_{1}+a_{2}]=[a_{1},a_{1}]+[a_{1},a_{2}]+[a_{2},a_{1}]+[a_{2},a_{2}]=\\
a_{2}+a_{3}\in\mathrm{\mathbf{Leib}}(L),
\end{gather*}
and hence, $\mathrm{\mathbf{Leib}}(L)=C_{2}$. Thus, we come to the following type of Leibniz algebras:
\begin{gather*}
L_{37}=Fa_{1}\oplus Fa_{2}\oplus Fa_{3}\ \mbox{where }[a_{1},a_{1}]=a_{2},[a_{1},a_{2}]=a_{3},\\
[a_{1},a_{3}]=[a_{2},a_{1}]=[a_{2},a_{2}]=[a_{2},a_{3}]=[a_{3},a_{1}]=[a_{3},a_{2}]=[a_{3},a_{3}]=0.
\end{gather*}
Note also that $\mathrm{\mathbf{Leib}}(L_{37})=[L_{37},L_{37}]=\zeta^{\mathrm{left}}(L_{37})=Fa_{2}\oplus Fa_{3}$, $\zeta(L_{37})=\zeta^{\mathrm{right}}(L_{37})=Fa_{3}$, $\mathrm{\mathbf{ncl}}(L_{37})=3$.

Suppose now that $\mathrm{\mathbf{ncl}}(L)=2$. Then, $L$ has an upper central series of a length 2:
$$\langle0\rangle=C_{0}\leqslant C_{1}\leqslant C_{2}=L.$$
Here, we have two possibilities: $\mathrm{\mathbf{dim}}_{F}(C_{1})=1$ or $\mathrm{\mathbf{dim}}_{F}(C_{1})=2$. Since $L/C_{1}$ is abelian, $\mathrm{\mathbf{Leib}}(L)\leqslant C_{1}$. Then, the fact that $\mathrm{\mathbf{dim}}_{F}(\mathrm{\mathbf{Leib}}(L))=2$ implies that $\mathrm{\mathbf{dim}}_{F}(C_{1})=2$.

Since $L$ is a not Lie algebra, there is an element $a_{1}$ such that $[a_{1},a_{1}]=a_{3}\neq0$. Then, $a_{1}\not\in\zeta(L)$, $a_{3}\in\mathrm{\mathbf{Leib}}(L)\leqslant\zeta(L)$. It follows that $[a_{1},a_{3}]=[a_{3},a_{1}]=[a_{3},a_{3}]=0$. We can see that $A=\langle a_{1}\rangle=Fa_{1}\oplus Fa_{3}$. In particular, $A\cap C_{1}=Fa_{3}$. We have $C_{1}=Fa_{2}\oplus Fa_{3}$ for some element $a_{2}$. Then, $L=A\oplus Fa_{2}$. Since $a_{2}\in\zeta(L)$, $Fa_{2}$ is an ideal of $L$. The choice of $a_{2}$ implies that $[a_{1},a_{2}]=[a_{2},a_{1}]=[a_{2},a_{3}]=[a_{3},a_{2}]=[a_{2},a_{2}]=0$. It follows that factor-algebra $L/Fa_{3}$ is abelian. Then, $\mathrm{\mathbf{Leib}}(L)\leqslant Fa_{3}$, in particular, $\mathrm{\mathbf{dim}}_{F}(\mathrm{\mathbf{Leib}}(L))=1$, and we obtain a contradiction. This contradiction shows that the case $\mathrm{\mathbf{ncl}}(L)=2$ is not possible.
\qed

\begin{thm}\label{T6}
Let $L$ be a non-nilpotent Leibniz algebra over a field $F$ having dimension $3$, which is not a Lie algebra. Suppose that $\zeta(L)\neq\langle0\rangle$ and that $\mathrm{\mathbf{dim}}_{F}(\mathrm{\mathbf{Leib}}(L))=2$. Then, $L$ is an algebra of the following type:
\begin{enumerate}
\item[\upshape] $\mathrm{\mathbf{Lei}}_{38}(3,F)=L_{38}$ is a cyclic Leibniz algebra, $L_{24}=Fa_{1}\oplus Fa_{2}\oplus Fa_{3}$ where $[a_{1},a_{1}]=a_{2}$, $[a_{1},a_{2}]=a_{2}+a_{3}$, $[a_{1},a_{3}]=[a_{2},a_{1}]=[a_{2},a_{2}]=[a_{2},a_{3}]=[a_{3},a_{1}]=[a_{3},a_{2}]=[a_{3},a_{3}]=0$, $\mathrm{\mathbf{Leib}}(L_{38})=\zeta^{\mathrm{left}}(L_{38})=[L_{38},L_{38}]=Fa_{2}\oplus Fa_{3}$, $\zeta^{\mathrm{right}}(L_{38})=\zeta(L_{38})=Fa_{3}$.
\end{enumerate}
\end{thm}
\pf We note that a Leibniz algebra of dimension 1 is abelian. Then, by our conditions, we obtain that $\mathrm{\mathbf{dim}}_{F}(\zeta(L))=1$. If we suppose that $\zeta(L)\cap\mathrm{\mathbf{Leib}}(L)=\langle0\rangle$, then $L=\zeta(L)\oplus\mathrm{\mathbf{Leib}}(L)$, so that $L$ is abelian, and we obtain a contradiction. This contradiction shows that $\zeta(L)\leqslant\mathrm{\mathbf{Leib}}(L)$. Put $C=\zeta(L)$ and let $c$ be an element such that $C=Fc$. Since $C\neq\mathrm{\mathbf{Leib}}(L)$, the factor-algebra $L/C$ is not a Lie algebra. Using the information above about the structure of the Leibniz algebras of dimension 2, we obtain that $L/C=F(b+C)\oplus F(d+C)$ where
\begin{gather*}
[d+C,d+C]=b+C,[d+C,b+C]=b+C,\\
[b+C,d+C]=[b+C,b+C]=C.
\end{gather*}
Without loss of generality, we may assume that $[d,d]=b$. Then, we have $[d,b]=b+\alpha c$, $[b,d]=0$ for some element $\alpha\in F$. The fact that $\mathrm{\mathbf{Leib}}(L)$ is abelian implies that $[b,b]=[b,c]=[c,b]=0$. Let $\lambda_{1}d+\lambda_{2}b+\lambda_{3}c$ be an arbitrary element of $L$. We have:
\begin{gather*}
[\lambda_{1}d+\lambda_{2}b+\lambda_{3}c,\lambda_{1}d+\lambda_{2}b+\lambda_{3}c]=\lambda_{1}^{2}b+\lambda_{1}\lambda_{2}(b+\alpha c)=\\
(\lambda_{1}^{2}+\lambda_{1}\lambda_{2})b+\lambda_{1}\lambda_{2}\alpha c.
\end{gather*}
Thus, we can see that if $\alpha=0$, then $\mathrm{\mathbf{Leib}}(L)=Fb$. In particular, we have $\mathrm{\mathbf{dim}}_{F}(\mathrm{\mathbf{Leib}}(L))=1$, and we obtain a contradiction. This contradiction shows that $\alpha\neq0$. Put $a_{1}=d$, $a_{2}=b$, $a_{3}=\alpha c$, then
$$[a_{1}-a_{2},a_{1}-a_{2}]=[a_{1},a_{1}]-[a_{2},a_{1}]-[a_{1},a_{2}]+[a_{2},a_{2}]=a_{2}-a_{2}-a_{3}=-a_{3}.$$
Thus, $\mathrm{\mathbf{Leib}}(L)=Fa_{2}\oplus Fa_{3}$, and we come to the following type of Leibniz algebras:
\begin{gather*}
L_{38}=Fa_{1}\oplus Fa_{2}\oplus Fa_{3}\ \mbox{where }[a_{1},a_{1}]=a_{2},[a_{1},a_{2}]=a_{2}+a_{3},\\
[a_{1},a_{3}]=[a_{2},a_{1}]=[a_{2},a_{2}]=[a_{2},a_{3}]=[a_{3},a_{1}]=[a_{3},a_{2}]=[a_{3},a_{3}]=0.
\end{gather*}
Note also that $\mathrm{\mathbf{Leib}}(L_{38})=\zeta^{\mathrm{left}}(L_{38})=[L_{38},L_{38}]=Fa_{2}\oplus Fa_{3}$, $\zeta(L_{38})=\zeta^{\mathrm{right}}(L_{38})=Fa_{3}$.
\qed

\begin{thm}\label{T7}
Let $L$ be a Leibniz algebra over a field $F$ having dimension $3$ and $L$ not be a Lie algebra. Suppose that $\zeta(L)=\langle0\rangle$ and that $\mathrm{\mathbf{dim}}_{F}(\mathrm{\mathbf{Leib}}(L))=2$. Then, $L$ is an algebra of one of the following types.
\begin{enumerate}
\item[\upshape(i)] $\mathrm{\mathbf{Lei}}_{39}(3,F)=L_{39}$ is a direct sum of ideal $B=Fa_{3}$ and cyclic non-nilpotent subalgebra $A=Fa_{1}\oplus Fa_{2}$, $[A,B]=Fa_{3}$, $[B,A]=\langle0\rangle$, so that $L_{39}=Fa_{1}\oplus Fa_{2}\oplus Fa_{3}$ where $[a_{1},a_{1}]=[a_{1},a_{2}]=a_{2}$, $[a_{1},a_{3}]=\beta a_{3}\ (\beta\neq0)$, $[a_{2},a_{1}]=[a_{2},a_{2}]=[a_{2},a_{3}]=[a_{3},a_{1}]=[a_{3},a_{2}]=[a_{3},a_{3}]=0$, $\mathrm{\mathbf{Leib}}(L_{39})=\zeta^{\mathrm{left}}(L_{39})=[L_{39},L_{39}]=Fa_{2}\oplus Fa_{3}$,\linebreak $\zeta(L_{39})=\zeta^{\mathrm{right}}(L_{39})=\langle0\rangle$.
\item[\upshape(ii)] $\mathrm{\mathbf{Lei}}_{40}(3,F)=L_{40}$ is a cyclic Leibniz algebra, so that $L_{40}=Fa_{1}\oplus Fa_{2}\oplus Fa_{3}$ where $[a_{1},a_{1}]=a_{2}$, $[a_{1},a_{2}]=a_{2}+\gamma a_{3}$ $(\gamma\neq0)$, $[a_{1},a_{3}]=\beta a_{3}\ (\beta\neq0)$, $[a_{2},a_{1}]=[a_{2},a_{2}]=[a_{2},a_{3}]=[a_{3},a_{1}]=[a_{3},a_{2}]=[a_{3},a_{3}]=0$, $\mathrm{\mathbf{Leib}}(L_{40})=\zeta^{\mathrm{left}}(L_{40})=[L_{40},L_{40}]=Fa_{2}\oplus Fa_{3}$, $\zeta^{\mathrm{right}}(L_{40})=\zeta(L_{40})=\langle0\rangle$.
\item[\upshape(iii)] $\mathrm{\mathbf{Lei}}_{41}(3,F)=L_{41}$ is a cyclic Leibniz algebra, so that $L_{41}=Fa_{1}\oplus Fa_{2}\oplus Fa_{3}$ where $[a_{1},a_{1}]=a_{2}$, $[a_{1},a_{2}]=\gamma a_{3}\ (\gamma\neq0)$, $[a_{1},a_{3}]=a_{3}$, $[a_{2},a_{1}]=[a_{2},a_{2}]=[a_{2},a_{3}]=[a_{3},a_{1}]=[a_{3},a_{2}]=[a_{3},a_{3}]=0$, $\mathrm{\mathbf{Leib}}(L_{41})=\zeta^{\mathrm{left}}(L_{41})=[L_{41},L_{41}]=Fa_{2}\oplus Fa_{3}$, $\zeta^{\mathrm{right}}(L_{41})=\zeta(L_{41})=\langle0\rangle$.
\item[\upshape(iv)] $\mathrm{\mathbf{Lei}}_{42}(3,F)=L_{42}$ is a cyclic Leibniz algebra, so that $L_{42}=Fa_{1}\oplus Fa_{2}\oplus Fa_{3}$ where $[a_{1},a_{1}]=a_{2}$, $[a_{1},a_{2}]=a_{3}$, $[a_{1},a_{3}]=\beta a_{2}+\gamma a_{3}$, $[a_{2},a_{1}]=[a_{2},a_{2}]=[a_{2},a_{3}]=[a_{3},a_{1}]=[a_{3},a_{2}]=[a_{3},a_{3}]=0$, $\mathrm{\mathbf{Leib}}(L_{42})=\zeta^{\mathrm{left}}(L_{42})=[L_{42},L_{42}]=Fa_{2}\oplus Fa_{3}$, $\zeta^{\mathrm{right}}(L_{42})=\zeta(L_{42})=\langle0\rangle$. Moreover, polynomial $X^{2}-\gamma X-\beta$ is irreducible over field $F$.
\end{enumerate}
\end{thm}
\pf Suppose first that $\mathrm{\mathbf{Leib}}(L)$ includes an ideal $K$ of dimension 1. Let $c$ be an element such that $K=Fc$. Since $K\neq\mathrm{\mathbf{Leib}}(L)$, the factor-algebra $L/K$ is not a Lie algebra. Using the above information about the structure of the Leibniz algebras of dimension 2, we obtain that $L/K=F(b+K)\oplus F(d+K)$ where
\begin{gather*}
[d+K,d+K]=b+K,[d+K,b+K]=b+K,\\
[b+K,d+K]=[b+K,b+K]=K
\end{gather*}
or
\begin{gather*}
[d+K,d+K]=b+K,\\
[d+K,b+K]=[b+K,d+K]=[b+K,b+K]=K.
\end{gather*}
Consider the first situation. Without loss of generality, we may assume that $[d,d]=b$. Then, we have $[b,d]=0$, $[d,b]=b+\alpha c$ for some element $\alpha\in F$. The fact that $\mathrm{\mathbf{Leib}}(L)$ is abelian implies that $[b,b]=[b,c]=[c,b]=0$. Since $\zeta(L)=\langle0\rangle$, $[d,c]=\beta c$ for some non-zero element $\beta\in F$. Put $a_{3}=\beta c$, then $K=Fa_{3}$ and $[b,a_{3}]=[a_{3},b]=[a_{3},d]=[a_{3},a_{3}]=0$, $[d,b]=b+\gamma a_{3}$ where $\gamma=\alpha\beta^{-1}$. Let $u=\lambda_{1}d+\lambda_{2}b+\lambda_{3}a_{3}$ be an arbitrary element of $L$. We have:
\begin{gather*}
[u,u]=[\lambda_{1}d+\lambda_{2}b+\lambda_{3}a_{3},\lambda_{1}d+\lambda_{2}b+\lambda_{3}a_{3}]=\\
\lambda_{1}^{2}[d,d]+\lambda_{1}\lambda_{2}[d,b]+\lambda_{1}\lambda_{3}[d,a_{3}]=\\
\lambda_{1}^{2}b+\lambda_{1}\lambda_{2}(b+\gamma a_{3})+\lambda_{1}\lambda_{3}\beta a_{3}=\\
(\lambda_{1}^{2}+\lambda_{1}\lambda_{2})b+(\gamma\lambda_{1}\lambda_{2}+\lambda_{1}\lambda_{3}\beta)a_{3}.
\end{gather*}
If we put $\lambda_{1}=1$, $\lambda_{2}=-1$, $\lambda_{3}=\beta^{-1}(1+\gamma)$, then we obtain $[u,u]=a_{3}$. Put $a_{1}=d$, $a_{2}=b$.

If $\alpha=0$ and hence, $\gamma=0$, then a subalgebra $Fa_{2}$ is an ideal of $L$ and $A=\langle a_{1}\rangle=Fa_{1}\oplus Fa_{2}$ is a cyclic subalgebra. Thus, we come to the following type of Leibniz algebras:
\begin{gather*}
L_{39}=Fa_{1}\oplus Fa_{2}\oplus Fa_{3}\ \mbox{where }[a_{1},a_{1}]=[a_{1},a_{2}]=a_{2},\\
[a_{1},a_{3}]=\beta a_{3}\ (\beta\neq0),\\
[a_{2},a_{1}]=[a_{2},a_{2}]=[a_{2},a_{3}]=[a_{3},a_{1}]=[a_{3},a_{2}]=[a_{3},a_{3}]=0.
\end{gather*}
Note also that $\mathrm{\mathbf{Leib}}(L_{39})=\zeta^{\mathrm{left}}(L_{39})=[L_{39},L_{39}]=Fa_{2}\oplus Fa_{3}$, $\zeta(L_{39})=\zeta^{\mathrm{right}}(L_{39})=\langle0\rangle$.

If $\alpha\neq0$, then $a_{3}=\gamma^{-1}([a_{1},a_{2}]-a_{2})$. It follows that $L$ is a cyclic algebra. Thus, we come to the following type of Leibniz algebras:
\begin{gather*}
L_{40}=Fa_{1}\oplus Fa_{2}\oplus Fa_{3}\ \mbox{where }[a_{1},a_{1}]=a_{2},\\
[a_{1},a_{2}]=a_{2}+\gamma a_{3}\ (\gamma\neq0),[a_{1},a_{3}]=\beta a_{3}\ (\beta\neq0),\\
[a_{2},a_{1}]=[a_{2},a_{2}]=[a_{2},a_{3}]=[a_{3},a_{1}]=[a_{3},a_{2}]=[a_{3},a_{3}]=0.
\end{gather*}
Note also that $\mathrm{\mathbf{Leib}}(L_{40})=\zeta^{\mathrm{left}}(L_{40})=[L_{40},L_{40}]=Fa_{2}\oplus Fa_{3}$, $\zeta(L_{40})=\zeta^{\mathrm{right}}(L_{40})=\langle0\rangle$.

Consider now a situation when $L/K=F(b+K)\oplus F(d+K)$ where
\begin{gather*}
[d+K,d+K]=b+K,\\
[d+K,b+K]=[b+K,d+K]=[b+K,b+K]=K.
\end{gather*}
Without loss of generality we may assume that $[d,d]=b$. Then we have $[b,d]=0$, $[d,b]=\alpha c$ for some element $\alpha\in F$. The fact that $\mathrm{\mathbf{Leib}}(L)$ is abelian implies that $[b,b]=[b,c]=[c,b]=0$. Since $\zeta(L)=\langle0\rangle$, $[d,c]=\beta c$ for some non-zero element $\beta\in F$. Put $a_{3}=c$, then $K=Fa_{3}$ and $[b,a_{3}]=[a_{3},b]=[a_{3},d]=[a_{3},a_{3}]=0$, $[d,b]=\alpha a_{3}$.

If $\alpha=0$, then $Fb$ lies in the center of $L$, and we obtain a contradiction.

Suppose that $\alpha\neq0$. Put $a_{1}=\beta^{-1}d$, then $[a_{1},c]=c$. Further $[a_{1},a_{1}]=[\beta^{-1}d,\beta^{-1}d]=\beta^{-2}b=a_{2}$. We have $[a_{2},a_{2}]=[a_{2},a_{1}]=[a_{2},c]=[c,a_{2}]=0$, $[a_{1},a_{2}]=[\beta^{-1}d,\beta^{-2}b]=\beta^{-3}\alpha c=\gamma c$. Thus we come to the following type of Leibniz algebra:
\begin{gather*}
L_{41}=Fa_{1}\oplus Fa_{2}\oplus Fa_{3}\ \mbox{where }[a_{1},a_{1}]=a_{2},\\
[a_{1},a_{2}]=\gamma a_{3}\ (\gamma\neq0),[a_{1},a_{3}]=a_{3},\\
[a_{2},a_{1}]=[a_{2},a_{2}]=[a_{2},a_{3}]=[a_{3},a_{1}]=[a_{3},a_{2}]=[a_{3},a_{3}]=0.
\end{gather*}
Note also that $\mathrm{\mathbf{Leib}}(L_{41})=\zeta^{\mathrm{left}}(L_{41})=[L_{41},L_{41}]=Fa_{2}\oplus Fa_{3}$, $\zeta(L_{41})=\zeta^{\mathrm{right}}(L_{41})=\langle0\rangle$.

Suppose now that $\mathrm{\mathbf{Leib}}(L)$ does not include proper non-zero ideals. Since $L$ is a non-Lie algebra, there is an element $d$ such that $[d,d]=b\neq0$. Then, $d\not\in\mathrm{\mathbf{Leib}}(L)$. Put $K=\mathrm{\mathbf{Leib}}(L)$. By our assumption, $[d,b]=c\not\in Fb$. The fact that $\mathrm{\mathbf{dim}}_{F}(K)=2$ implies that $K=Fb\oplus Fc$. Then, $[d,c]=\beta b+\gamma c$ for some elements $\beta,\gamma\in F$. The mapping $x\rightarrow[d,x]$, $x\in K$ is linear. Our conditions imply that a polynomial $X^{2}-\gamma X-\beta$ is irreducible over a field $F$. Put $a_{1}=d$, $a_{2}=b$, $a_{3}=c$. Thus, we come to the following type of Leibniz algebras:
\begin{gather*}
L_{42}=Fa_{1}\oplus Fa_{2}\oplus Fa_{3}\ \mbox{where }[a_{1},a_{1}]=a_{2},[a_{1},a_{2}]=a_{3},\\
[a_{1},a_{3}]=\beta a_{2}+\gamma a_{3},\\
[a_{2},a_{1}]=[a_{2},a_{2}]=[a_{2},a_{3}]=[a_{3},a_{1}]=[a_{3},a_{2}]=[a_{3},a_{3}]=0.
\end{gather*}
Note also that $\mathrm{\mathbf{Leib}}(L_{42})=\zeta^{\mathrm{left}}(L_{42})=[L_{42},L_{42}]=Fa_{2}\oplus Fa_{3}$, $\zeta(L_{42})=\zeta^{\mathrm{right}}(L_{42})=\langle0\rangle$, polynomial $X^{2}-\gamma X-\beta$ is irreducible over field $F$.
\qed

\end{document}